\documentclass[a4paper]{elsarticlemod}

\usepackage[dvipsnames]{xcolor}
\usepackage{graphicx}
\usepackage{amssymb}
\usepackage{amsmath}
\usepackage{amsthm}
\usepackage{bm}
\usepackage{subfigure}
\usepackage{caption} 
\DeclareCaptionFormat{cont}{#1 (cont.)#2#3\par} 
\allowdisplaybreaks 
\usepackage[T1]{fontenc}
\usepackage{algorithm2e}
\usepackage{soul}
\usepackage{stmaryrd}
\SetSymbolFont{stmry}{bold}{U}{stmry}{m}{n}
\usepackage{enumerate}
\usepackage{stackengine}
\usepackage[colorinlistoftodos, obeyFinal]{todonotes}
\usepackage{multirow}
\usepackage{hhline}
\usepackage{textcomp}
\usepackage{MnSymbol}
\usepackage{colortbl}
\usepackage{tikz}
\usetikzlibrary{external,shapes.misc}
\usepackage{pgfplots}
\pgfplotsset{compat=1.7}
\usepackage{etoolbox}
\usetikzlibrary{pgfplots.patchplots}
\usetikzlibrary{pgfplots.colormaps}
\usetikzlibrary{pgfplots.groupplots}
\usetikzlibrary{positioning,arrows,mindmap,backgrounds,shapes,decorations.fractals,calc,intersections,decorations.pathmorphing,patterns}

\pgfplotsset{stress plot style/.style={
		width=0.48\textwidth,
		y label style={at={(axis description cs:-0.10,.5)},anchor=south},
		ylabel={$x_3$ [mm]}}}

\newif\ifrecompiletikz
\recompiletikzfalse

\ifrecompiletikz
\usetikzlibrary{external}
\tikzexternaldisable
\else
  \tikzexternaldisable
\fi

\tikzset{cross/.style={cross out, minimum size=2*(#1-\pgflinewidth), inner sep=0pt, outer sep=0pt},
cross/.default={1pt}}

\newcommand{\blueline}{\textcolor{blue}{$\boldsymbol{\leftrightline}$}}
\newcommand{\redcrosses}{\textcolor{red}{$\boldsymbol{\times}$}}
\DeclareRobustCommand{\blackcircles}{\protect\raisebox{0.0cm}{\tikz{\node[draw=black,line width=0.8pt,scale=0.5,circle,fill=white](){};}}}
\DeclareRobustCommand{\blackcirclesfull}{\protect\raisebox{0.0cm}{\tikz{\node[draw=black,line width=0.8pt,scale=0.3,circle,fill=black](){};}}}

\DeclareRobustCommand{\bluesolidcircle}{\protect\raisebox{0.0cm}{\tikz{\node[draw=blue,line width=0.8pt,scale=0.6,circle,fill=none](){}; \draw[draw=blue,line width=0.8pt] (-0.25,0)--(0.25,0);}}}
\DeclareRobustCommand{\blacksolidx}{\protect\raisebox{0.0cm}{\tikz{\node[draw=black,line width=0.8pt,scale=0.6,cross=4pt,fill=none](){}; \draw[draw=black,line width=0.8pt] (-0.25,0)--(0.25,0);}}}
\DeclareRobustCommand{\redsolidcross}{\protect\raisebox{0.0cm}{\tikz{\node[draw=red,line width=0.8pt,scale=0.6,cross=3.5pt,rotate=45,fill=none](){}; \draw[draw=red,line width=0.8pt] (-0.25,0)--(0.25,0);}}}

\pgfplotsset{select coords between index/.style 2 args={
		x filter/.code={
			\ifnum\coordindex<#1\fi
			\ifnum\coordindex>#2\fi
		}
}}

\newtheorem{remark}{Remark}

\tikzset{
	dot/.style = {circle, fill=blue, minimum size=#1,
		inner sep=0pt, outer sep=0pt},
	dot/.default = 2pt 
}
\tikzset{
	dot2/.style = {circle, fill=black, minimum size=#1,
		inner sep=0pt, outer sep=0pt},
	dot2/.default = 2pt 
}

\definecolor{majorelleblue}{rgb}{0.38, 0.31, 0.86}
\definecolor{frenchblue}{rgb}{0.0, 0.45, 0.73}
\definecolor{deepchestnut}{rgb}{0.73, 0.31, 0.28}
\definecolor{dgreen}{RGB}{0,139,0}
\definecolor{forestgreen}{rgb}{0.0,0.27,0.13}
\definecolor{darkviolet}{rgb}{0.58, 0.0, 0.83}
\definecolor{coolgrey}{rgb}{0.55, 0.57, 0.67}
	
\newcommand{\diff}{\textrm{\normalfont d}}

\newcommand{\strain}{\bm{\varepsilon}}
\newcommand{\stress}{\bm{\sigma}}
\newcommand{\CC}{\mathbb{C}}
\newcommand{\CCo}{\overline{\mathbb{C}}}
\newcommand{\RR}{\bm{\mathsf{R}}}
\newcommand{\uu}{\bm{\mathsf{u}}}

\newcommand{\one}{\bm{I}}

\newcommand{\T}[1]{{\bm{#1}}}

\def\e{\bm{e}}
\def\a{\bm{a}}

\def\E{\bm{E}}
\def\g{\bm{g}}
\newcommand{\tensorm}{%
	\setstackgap{S}{0.4ex}%
	\mathrel{\Shortstack{{.} {.} {.}}}}

\makeatletter
\newcommand{\vast}{\bBigg@{17.0}}
\newcommand{\thickhline}{
	\noalign {\ifnum 0=`}\fi \hrule height 1pt
	\futurelet \reserved@a \@xhline
}
\newcolumntype{?}{!{\vrule width 1pt}}
\makeatother
\newcommand\Tstrut{\rule{0pt}{2.6ex}}         
\newcommand\Bstrut{\rule[-0.9ex]{0pt}{0pt}}   

\newcommand{\ie}{{i.e.,~}}
\newcommand{\eg}{{e.g.,~}}

\def\rmd{{\mathrm{d}}}

\def\be{\begin{equation}}
\def\ee{\end{equation}}
\def\ba{\begin{array}}
\def\ea{\end{array}}

\begin{document}	
	
	\begin{frontmatter}
		
        \title{Efficient equilibrium-based stress recovery for isogeometric laminated curved structures}
		\author[pavia]{Alessia Patton\corref{cor1}}
		\author[losanne]{Pablo Antol\'in}
		\author[munchen]{Josef Kiendl}
		\author[pavia]{Alessandro Reali}
		\address[pavia]{Department of Civil Engineering and Architecture - University of Pavia\\
			via Ferrata 3, 27100, Pavia, Italy}  
		\address[losanne]{Institute of Mathematics - \'Ecole Polytechnique F\'ed\'erale de Lausanne\\
			CH-1015 Lausanne, Switzerland} 
		\address[munchen]{Department of Civil Engineering and Environmental Sciences - Universit\"at der Bundeswehr M\"unchen\\Werner-Heisenberg-Weg 39, 85577 Neubiberg, Germany}
		\cortext[cor1]{\\Corresponding author. Email: alessia.patton01@universitadipavia.it}
			
	\begin{abstract}
	  This work focuses on an efficient stress recovery procedure for laminated composite curved structures, which relies on Isogeometric Analysis (IGA) and equilibrium. Using a single
	  element through the thickness in combination with a calibrated layerwise integration rule or a homogenized approach, the 3D solid isogeometric modeling grants an inexpensive and accurate approximation in terms of displacements (and their derivatives) and in-plane stresses, while through-the-thickness stress components are poorly approximated. Applying a further post-processing step, an accurate out-of-plane stress state is also recovered, even from a coarse displacement solution. This is based on a direct integration of the equilibrium equations in strong form, involving high order derivatives of the displacement field. Such a continuity requirement is fully granted by IGA shape function properties. The post-processing step is locally applied, which grants that no additional coupled terms appear in the equilibrium, allowing for a direct reconstruction without the need to further iterate to resolve the out-of-balance momentum equation. Several numerical results show the good performance of this approach particularly for composite stacks with significant radius-to-thickness ratio and number of plies. In particular, in the latter case, where a layerwise technique employing a number of degrees of freedom directly proportional to the number of plies would be much more computationally demanding, the proposed method can be regarded as a very appealing alternative choice.
	\end{abstract}

\end{frontmatter}
    
	\section{Introduction}\label{sec:introduction}
    Composite structures are formed by a combination of two or more components such that enhanced mechanical properties can be obtained \cite{Gibson1994,Jones1998,Reddy2004}. Being characterized by a high stiffness and strength despite maintaining a reduced weight, the interest for composites in the engineering community has constantly grown in recent years, especially in the aerospace and automotive industries.

In this work we focus on laminated structures, which are made of a collection of laminae that have a preferred high-strength direction. The layers are stacked and subsequently glued together such that the orientation of the high-strength direction varies with each successive layer giving designers the flexibility to tailor laminate properties. Nevertheless, it is a well known fact that composites are prone to damage under even simple loading conditions due to comparatively poor strength in the out-of-plane direction. One of the most common failure modes of composite structures is referred to as delamination, which consists of an interface crack between two adjacent plies \cite{Sridharan2008}. The most common sources of delamination are material and structural discontinuities that give rise to relevant interlaminar stresses \cite{Mittelstedt2007}. The interlaminar stress level is strongly dependent on the composite stacking sequence \cite{Pagano1971,Pagano1973} and the mismatch of engineering properties between adjacent plies \cite{Jones1998}. Therefore, the capability of accurately predicting out-of-plane stresses is of paramount importance.

Over the past decades, various theories have been proposed to model laminated composite plates and shells. These can be classified into three major categories: Three dimensional (3D) elasticity, equivalent single-layer (ESL), and Layerwise (LW) theories \cite{Li2020,Liew2019}.
Considering the laminated composite structure as a solid without any special treatment of the stacking sequence, three dimensional elasticity theories, despite being accurate, lead to heavy computations. 
In an attempt to reduce the computational effort, ESL treat instead the 3D laminate as an equivalent single ply, such that the number of unknowns is independent of the number
of layers. The major drawback of ESL derives from the continuous assumption of the displacement field approximation, predicting continuous transverse strains, which, together with ply-wise different material properties, leads to discontinuous interlaminar stresses.
Nonetheless, ESL can provide good results
for the global responses of thin laminated composite
plates and shells, while they fail to capture local responses, such as the distribution of interlaminar stresses, in particular in the case of thick laminates.
Numerous theories based on the ESL concept have
been proposed \cite{Khandan2012,Kreja2011,Reddy2004}. While \textit{classical laminate theory} neglects interlaminar strains, the \textit{first order shear deformation theory} (FSDT) considers the transverse shear strains to be constant through-the-thickness and therefore requires correction factors, which are difficult to determine for arbitrarily laminated structures. In an attempt to address this issue, higher-order ESL theories have been proposed \cite{Lo1977,Reddy1984}, which rely on higher-order polynomials in the expansion of the displacement components through the thickness of the laminate. Nevertheless, despite the fact that they do not need correction factors, higher-order theories introduce additional unknowns that often do not have a clear physical meaning and increase the overall computational demand \cite{Reddy2004}.

To overcome the limitations of ESL and 3D elasticity models, LW theories can be considered as a viable option. Those approaches further distinguish between displacement-based strategies, for which the displacement field is the primal variable, and mixed LW theories, that comprise both displacements and transverse stresses as unknowns. Displacement-based LW theories consider an independent displacement field in every single ply and impose $C^0$-continuity at the layer interfaces, decreasing the number of unknowns \cite{Liu1996,Reddy2004}. Therefore, LW theories naturally fulfil the requirements which grant a good approximation of the out-of-plane response through the thickness directly from the constitutive equations. On the other hand, mixed LW models satisfy the interlaminar continuity of transverse stresses a priori. The main limitation of the LW theories is that the number of degrees of freedom is directly proportional to the number of layers, unavoidably leading to high computational costs. Further approximation theories comprise the ``Carrera Unified Formulation'', which allows to directly derive a large number of higher-order laminate models with a variable number of unknowns referring to few fundamental nuclei \cite{Carrera2011}. Recent works employ Isogeometric Analysis (IGA) \cite{Hughes2005} in the context of LW theories, benefiting from IGA predictive capabilities while keeping a significant low computational cost \cite{Dufour2018,Guo2014,Guo2015,Thai2016}.

IGA aims at bridging the gap between Finite Elements Analysis and Computer Aided Design (CAD), exploiting simulation methods that use the same basis functions both for the representation of the geometric computational domains and for the numerical simulations. This is accomplished by employing CAD basis functions (like B-Splines and NURBS), which present higher continuity properties than classical finite alement shape functions, allowing for an exact representation of complex, curved geometries, while typically leading to improved convergence and significantly faster overall simulations. IGA proved to be successful in a wide variety of solid and structural problems including - as pioneer works in the field - structural vibrations \cite{Cottrell2006}, nearly incompressible linear and non-linear elasticity and plasticity \cite{Elguedj2008}, phase-field modeling \cite{Borden2012}, plate \cite{Beiraodaveiga2012} and shell structures \cite{Kiendl2009,Benson2010}. Recent contributions specifically in the context of shell modeling comprise formulations allowing for large-strain plastic deformation \cite{Ambati2018}, the analysis of geometrically nonlinear elastic shells \cite{Leonetti2019}, and novel approaches for alleviating shear and membrane locking phenomena in solid shells \cite{Antolin2020b}. In the context of composite modeling, 2D IGA finite element formulations have been proposed in \cite{Kapoor2013,Nguyen-Xuan2013},
with some of them employing enhanced shell and plate theories \cite{Remmers2015,Adams2020}. Other approaches use a combination of IGA and
the Carrera Unified Formulation as in \cite{Alesadi2018}. IGA higher order shear deformation theories have been instead proposed for composite beams \cite{Shafei2020}, plates \cite{Thai2015,Thanh2019TWS}, and shells \cite{Casanova2013,Tornabene2017,Faroughi2020}, as well as for several examples of functionally graded (FG) plates (see, \eg \cite{Thanh2019CS,PhungVan2017} and references therein). With regard to IGA FG shell analysis,
investigations are restricted to a few studies including first \cite{Nguyen2018} and high order shear deformation theories \cite{VanDo2020}.

In this paper, we focus on interlaminar stress modeling of 3D curved laminated composites in the framework of higher order and higher continuity NURBS, extending the procedure presented in \cite{Dufour2018,Patton2019,Patton2020} for solid plates, which proved to be effective also in the context of Radial Basis Functions \cite{Chiappa2020}. The effectiveness of the isogeometric paradigm in the
modeling of composite laminates is demonstrated through several numerical examples which feature a 3D cylindrical shell under bending. The proposed strategy can be regarded as a two-step procedure: First, the structure is modeled using only one element through the thickness and a highly continuous displacement field is obtained via (coarse) solid isogeometric computations, which rely either on a layerwise integration rule or a homogenized approach; then, the interlaminar stresses are directly  recovered imposing equilibrium, starting from the same displacement solution and computing the necessary high-order derivatives.

Stress recovery theory examples in the context of finite element plates can be found in original works such as \cite{Pryor1971}, which proposes high order elements, still advocating for a generalization of the process that would include a separate angle of rotation for each layer in the laminate. Moreover, an average rotation of the mid-plane through the entire thickness of the plate is assumed, which requires shear coefficients depending on the section shape and also the recovery of the out-of-plane normal stress component is not discussed. In \cite{Engblom1985} instead, the proposed element formulation comprises different interpolation for different unknowns, resulting in a-not-so-straightforward procedure as far as implementation. Furthermore, this work emphasizes that a reduced numerical integration is strongly recommended especially in the case of thin geometries. In practice though, only the transverse shear components are integrated with reduced order, operation which may affect the physical behaviour of the element by introducing spurious zero energy modes. Another approach \cite{Ubertini2004}, instead, recovers stresses by minimizing the complementary energy functional associated to a separate patch system (on which finite element displacements are prescribed along the boundaries), over a set of stress fields which satisfy a priori interior equilibrium within the patch. This hybrid stress approach, tested for one- and two-dimensional isotropic elasticity problems, has been further analyzed in the context of FSDT laminated plates in \cite{Daghia2008}.
	
The herein introduced equilibrium-based technique, taking its inspiration from the aforementioned works, is more straightforward and extends the post-processing to more complex geometries, relying on the following peculiarities:
\begin{itemize}
	\item [-] An accurate in-plane solution can be obtained using only one element through the thickness.
	\item [-] The stress recovery is based on an accurate evaluation of higher-order in-plane stress derivatives which in turn require a higher displacement regularity, which is fully granted by IGA properties.
	\item [-] While the displacement field is obtained in a global framework, the interlaminar stresses are recovered locally, which grants that no further coupled terms arise in the equilibrium equations. This allows to handle the increasing geometry complexity of the curved solid case without the need to iterate solving the out-of-balance system. 
	\item [-] The ability to compute a valid guess for interlaminar shear stress boundary conditions allows to successfully recover also the normal out-of-plane stress component.
	\item [-] The possibility to exploit IGA geometry modeling features (\eg conic sections can be exactly constructed).
\end{itemize}

The structure of the paper is as follows. The governing equations for the 3D orthotropic elastic case are outlined in Section~\ref{sec:governing_eqq}, after introducing a global and a local description of the needed kinematics and constituive relations. Then, in Section~\ref{sec:recovery} the proposed post-processing strategy based on equilibrium is detailed starting from the approximated global displacement field, which is obtained using a 3D IGA single-element approach and a calibrated layerwise integration rule or a homogenized approach, introduced in Section~\ref{sec:numerical_strategies}. Several benchmarks comprising different mean radius-to-thickness ratios and numbers of layers are considered to showcase the effectiveness of the method in Section~\ref{sec:numerical_tests}. We finally draw our conclusions in Section~\ref{sec:conclusions}.
	
	\section{Governing equations for the orthotropic elastic case}\label{sec:governing_eqq}
	In this section, we study the equations that govern the 3D orthotropic elasticity problem. For this purpose, we introduce a global and a local cartesian reference system, as well as the associated kinematic quantities and constitutive relations needed to address the curved geometry problem.
The reason why we underline this global-local dual description of engineering quantities lies in our double-step stress recovery approach: First we tackle the problem looking for a global displacement solution which is then post-processed locally to reconstruct a proper out-of plane stress field from the in-plane one.
Finally, we detail the governing equations for the orthotropic solid case both at the strong and weak form level.

\subsection{Kinematics: A global and local perspective}\label{subsec:kin}
Let us consider an open bounded domain $\Omega\subset\mathbb{R}^3$ representing an elastic three-dimensional body, defined as a spline parameterization, such that $\Omega$ is the image of a unit cube $\hat\Omega = [0,1]^3$ through the map $\bm{F}:(\xi^1,\xi^2,\xi^3) \in \hat\Omega\to\bm{X}\in\Omega$, \ie $\Omega = \bm{F}(\hat\Omega)$.
\begin{figure}[htpb!]
	\centering
	\includegraphics[width=0.75\textwidth]{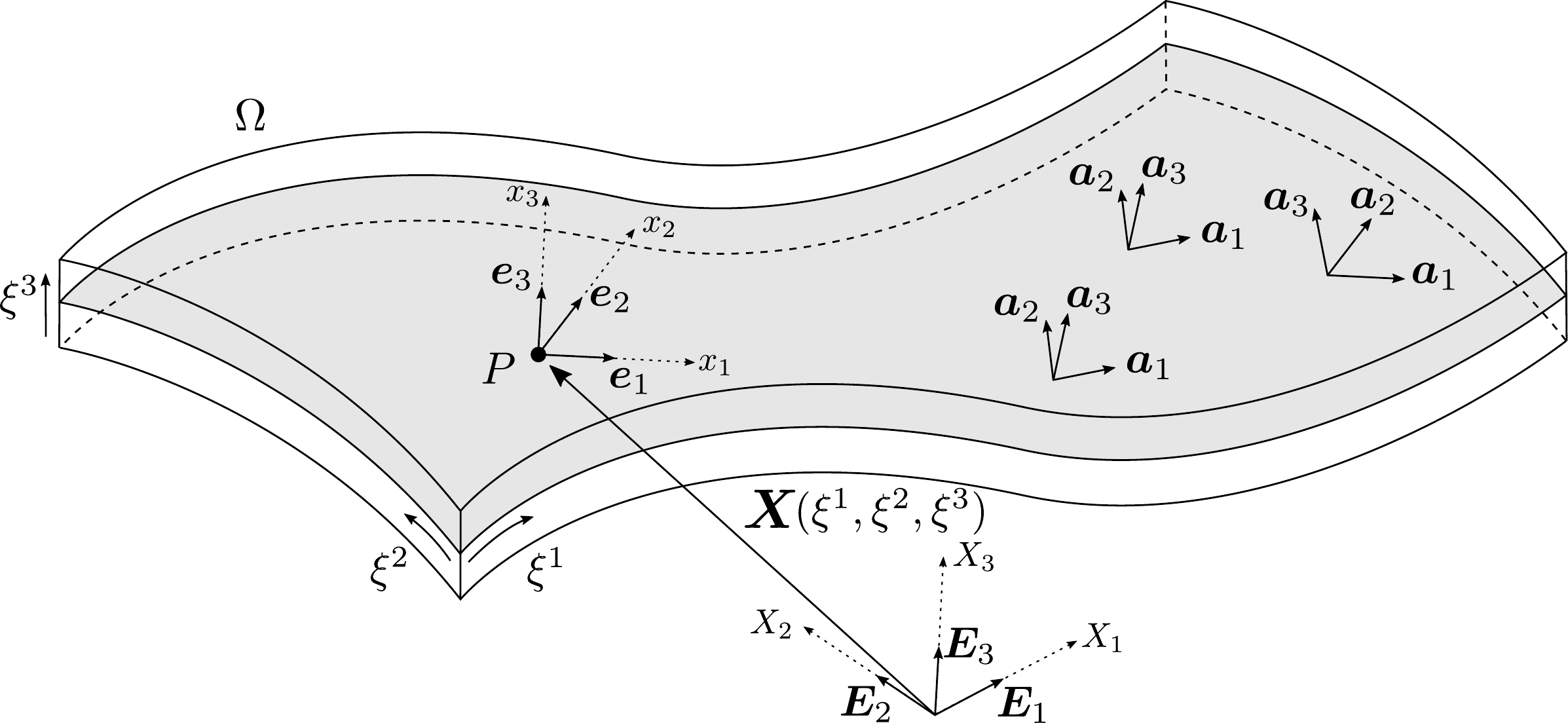}
	\caption{Global and local reference cartesian systems associated to the curved structure.}
	\label{fig:001}
\end{figure}

Let us also introduce a global reference cartesian system $\{X_1,X_2,X_3\}$ and its associated orthonormal basis $\{\E_1,\E_2,\E_3\}$,
as well as a local orthornormal basis $\{\a_1,\a_2,\a_3\}$, defined for every point of the body (see Figure~\ref{fig:001}).
This local basis is described according to the structure's mid surface (represented in gray color in Figure~\ref{fig:001}), such as vectors $\a_1$ and $\a_2$ define the in-plane directions of the structure and $\a_3$ the out-of-plane direction.
Thus, the same local basis $\{\a_1,\a_2,\a_3\}$ is associated to all the points along the thickness that have the same in-plane parametric coordinates $(\xi^1,\xi^2)$.
Apart from these requirements, the basis $\{\a_1,\a_2,\a_3\}$ can be chosen freely, nevertheless, in Appendix \ref{subsec:localBasis_der} we propose a particular choice based on an orthonormalization of the covariant basis of $\bm{F}$.

As already introduced in Section~\ref{sec:introduction}, the out-of-plane stress components will be recovered starting from the in-plane ones at selected points.
Once a point is set, and therefore defined by its in-plane coordinates $(\xi^1,\xi^2)$, we post-process a precise stress description along the third coordinate $\xi^3$.
Thus, for a given point of interest $P$ we introduce a fixed local cartesian reference system $\{x_1,x_2,x_3\}$ and its associated orthonormal basis $\{\e_1,\e_2,\e_3\}$ (see Figure~\ref{fig:001}).
This fixed local system is an auxiliary reference system that will be used during the calculations involved in the recovery process of the stress at $P$.
The basis $\{\e_1,\e_2,\e_3\}$ is just a static snapshot of the moving basis $\{\a_1,\a_2,\a_3\}$ evaluated at $P$.
Thus, due to the fact that $\{\e_1,\e_2,\e_3\}$ is fixed, its derivatives with respect to both reference systems vanish:
\begin{equation} \label{eq:fixed_local_basis}
    \frac{\diff \e_\alpha}{\diff X_i} = \bm{0}\,,\quad \frac{\diff \e_\alpha}{\diff x_\beta}= \bm{0}\,
\end{equation}
and the same applies to its subsequent derivatives. Nevertheless, this is not true, in general, for the moving local basis $\{\a_1,\a_2,\a_3\}$.

In the expression above, and hereinafter, latin indices are used for quantities expressed in the basis $\{\E_1,\E_2,\E_3\}$, whereas the greek indices refer to quantities expressed in the local bases $\{\e_1,\e_2,\e_3\}$ and $\{\a_1,\a_2,\a_3\}$, or parametric coordinates.
In addition, all indices span from 1 to 3 and we adopt Einstein's notation, \ie repeated indices imply the summation of the involved components unless otherwise stated.
For the sake of conciseness, henceforward, the global and local reference systems and bases are denoted, respectively, as $\{X_i\}$, $\{x_\alpha\}$, $\{\E_i\}$, $\{\e_\alpha\}$, and $\{\a_\alpha\}$.

The global basis $\{\E_i\}$ can be expressed in terms of the basis $\{\e_\alpha\}$, and viceversa, as:
\begin{subequations} \label{eq:bases_dot} \begin{align}
   \E_i &= C_{i\alpha}\,\e_\alpha\,,\\
   \e_\alpha &= C_{i\alpha}\,\E_i\,,
\end{align} \end{subequations}
where the basis change operator $C_{i\alpha}$ is defined as:
\begin{align} \label{eq:metric}
   C_{i\alpha} = \frac{\diff X_i}{\diff x_\alpha} = \frac{\diff x_\alpha}{\diff X_i} = \E_i \cdot \e_\alpha\,.
\end{align}
Due to the fact that $\{\E_i\}$ and $\{\e_\alpha\}$ are orthonormal bases, the change of basis operator is orthogonal.
It is worth noting that, due to the fact that both reference systems $\{X_i\}$ and $\{x_\alpha\}$ are fixed, $C_{i\alpha}$ is constant, i.e.:
\begin{align}\label{eq:metric_ders}
    \frac{\diff C_{i\alpha}}{\diff X_j} = 0\,, \quad
    \frac{\diff C_{i\alpha}}{\diff x_\beta} = 0\,,
\end{align}
and the same applies to subsequent derivatives.

Likewise, the global basis $\{\E_i\}$ can be expressed in terms of the local moving basis $\{\a_\alpha\}$ as:
\begin{align}
   \E_i &= D_{i\alpha}\,\a_\alpha\,,
\end{align}
where the basis change operator $D_{i\alpha}$ is defined as:
\begin{align} \label{eq:metric_2}
   D_{i\alpha} = \E_i \cdot \a_\alpha\,.
\end{align}
In general, $D_{i\alpha}$ is not constant, and therefore its derivatives with respect to $\{X_i\}$ or $\{x_\alpha\}$ do not vanish.

\subsection{Constitutive relations}\label{subsec:const_eq}
Considering the global basis $\{\E_i\}$ introduced above, we can express the 3D elastic displacement field as $\bm{u}=\tilde u_i \, \bm{E}_i$, 
where with $(\tilde{.})$ we refer to quantities expressed in the global basis $\{\E_i\}$.
We shall now assume small strains and small displacements, such that strains are given by the tensor $\strain = \T{\nabla}^s\bm{u}$, whose components in the global basis are
$\strain = \tilde \varepsilon_{ij} \, \bm{E}_i\otimes\bm{E}_j$,  where:
\begin{align}\label{eq:strain_components}
\tilde \varepsilon_{ij} = \frac{1}{2}\left( \frac{\diff \tilde u_i}{\diff X_j} + \frac{\diff \tilde u_j}{\diff X_i} \right)\,.
\end{align}

The displacement and strain fields can be expressed as well in the local basis $\{\e_\alpha\}$ as $\bm{u} = u_\alpha \e_\alpha$ and $\strain = \varepsilon_{\alpha\beta} \e_\alpha\otimes\e_\beta$, respectively.
Applying the change of basis defined in Equation~\eqref{eq:bases_dot}, the components $u_\alpha$ and $\varepsilon_{\alpha\beta}$ are computed as:
\begin{align}
u_\alpha &= \tilde u_i\, C_{i\alpha}\,, \\
\varepsilon_{\alpha\beta} &= \tilde{\varepsilon}_{ij}\, C_{i\alpha}\, C_{j\beta}\,.
\end{align}

In the case of linear elasticity, we introduce the stress field as
\begin{equation}\label{eq:stress_strain}
\stress = \CC : \strain\,,
\end{equation}
where $\CC$ is the fourth order material elasticity tensor, that can be expressed at every point using the global $\{\E_i\}$ or local $\{\a_\alpha\}$ bases
as:
\begin{align}\label{eq:CC_local}
\CC = \tilde{\CC}_{ijkl}\,\E_i\otimes\E_j\otimes\E_k\otimes\E_l = \CC_{\alpha\beta\gamma\delta}\,\a_\alpha\otimes\a_\beta\otimes\a_\gamma\otimes\a_\delta\,,
\end{align}
whose components are related through the change of basis operator~\eqref{eq:metric_2} as:
\begin{align} \label{eq:CC_change_basis}
\tilde\CC_{ijkl}=D_{i\alpha}\,D_{j\beta}\,D_{k\gamma}\,D_{l\delta}\,\CC_{\alpha\beta\gamma\delta}\,.
\end{align}

Focusing on the linear orthotropic elasticity case, $\CC$ may be expressed in a more convenient way using the local basis $\{\a_i\}$.
In fact, in the case that ply fiber directions locally coincide with the basis $\{\a_i\}$, the number of elastic coefficients of $\CC_{\alpha\beta\gamma\delta}$ reduces to nine. In Voigt notation, $\CC$ can be expressed in terms of engineering constants as
\begin{equation}\label{eq:CC_coeffs}
{\CC=\renewcommand\arraystretch{1.7}\vast{[}\begin{matrix}\CC_{11} & \CC_{12} & \CC_{13} & 0 & 0 & 0\\
	& \CC_{22} & \CC_{23} & 0 & 0 & 0\\
	&  & \CC_{33} & 0 & 0 & 0\\
	&  \text{symm.}&  & \CC_{44} & 0 & 0\\
	&  &  &  & \CC_{55} & 0\\
	&  &  &  &  & \CC_{66}\\
	\end{matrix}\vast{]}=\vast{[}\begin{matrix}
	\cfrac{1}{E_{1}} & -\cfrac{\nu_{12}}{E_{1}} & -\cfrac{\nu_{13}}{E_{1}} & 0 & 0 & 0\\
	& \cfrac{1}{E_{2}} & -\cfrac{\nu_{23}}{E_{2}} & 0 & 0 & 0\\
	&  & \cfrac{1}{E_{3}} & 0 & 0 & 0\\
	&  \text{symm.}&  & \cfrac{1}{G_{23}} & 0 & 0\\
	&  &  &  & \cfrac{1}{G_{13}} & 0\\
	&  &  &  &  & \cfrac{1}{G_{12}}\\
	\end{matrix}\vast{]}^{-1}\,.}
\end{equation} 

Like the displacement or the strain, also the stress can be expressed either in the global basis as $\stress=\tilde{\sigma}_{ij}\,\E_i\otimes\E_j$ or the local one as
$\stress=\sigma_{\alpha\beta}\,\e_\alpha\otimes\e_\beta$, where:
\begin{subequations}\label{eq:stress_local_global_coefs}
\begin{align}
\tilde\sigma_{ij}&=\tilde{\CC}_{ijkl}\,\tilde\varepsilon_{kl}\,,\label{eq:stress_global_coefs}\\
\sigma_{\alpha\beta}&=\tilde{\sigma}_{ij}C_{i\alpha}\,C_{j\beta}\,.\label{eq:stress_local_coefs}
\end{align}
\end{subequations}
	
\subsection{Strong form}\label{subsec:strong_form}
The considered elastic three-dimensional body is subjected to volume forces $\T{b}$, prescribed displacements $\bm{h}$ on
the Dirichlet portion $\Gamma_{D}$ of the boundary, and prescribed tractions $\boldsymbol{t}$ acting on the remaining
Neumann portion $\Gamma_{N}$, such that $\Gamma_{N}\cup\Gamma_{D}=\partial\Omega$ and $\Gamma_{N}\cap\Gamma_{D}=\emptyset$.

Hereby we recall the momentum balance equation in strong form and the corresponding boundary conditions for the linear elasticity problem 
\begin{subequations}\label{eq:strong_form}
	\begin{alignat}{3}
	& \T{\nabla} \cdot \T{\sigma} + \T{b} = \T{0}                                   \quad\quad &    \quad & \quad\hbox{in}\quad\Omega\,,
	\label{eq:equil}\\
	& \stress \cdot \boldsymbol{n}=\boldsymbol{t}   \quad\quad &    \quad & \quad\hbox{on}\quad\Gamma_N\,,\label{eq:nbc}\\
	& \boldsymbol{u} = \bm{h}                                                   \quad\quad &    \quad & \quad\hbox{on}\quad\Gamma_D\,,\label{eq:dbc}
	\end{alignat}
\end{subequations}
 where $\T{\nabla}\cdot$ represents the divergence operator computed with respect to the global cartesian reference system $\{X_i\}$ and $\T{n}$ is the outward normal unit vector. The stress tensor $\stress = \tilde{\sigma}_{ij}\, \E_i\otimes\E_j$ has its components defined as in Equation~\eqref{eq:stress_global_coefs}, while $\boldsymbol{b}=\tilde{b}_i\,\E_i$, $\boldsymbol{t}=\tilde{t}_i\,\E_i$, and $\bm{h}=\tilde{h}_i\,\E_i$, with $\tilde{b}_i$, $\tilde{t}_i$, and $\tilde{h}_i$, being, respectively, the body force, traction, and imposed displacement components in the basis $\{\E_i\}$.
 
Considering the system of equations~\eqref{eq:strong_form}, the term $\T{\nabla} \cdot \T{\sigma}$ needs to be further detailed in terms of global reference system $\{X_i\}$. Therefore, we need to develop equation~\eqref{eq:equil} expressing the stress tensor in terms of the constitutive relation in Equation~\eqref{eq:stress_strain} 
\begin{equation}\label{eq:stress_div}
\T{\nabla} \cdot \T{\sigma} = \T{\nabla} \cdot \CC : \strain + \CC \tensorm \T{\nabla}\strain\,, 
\end{equation}
where $\tensorm$ represents the triple contraction operator and
\begin{equation}\label{eq:stress_div_global}
\T{\nabla} \cdot \T{\sigma} = (\tilde{\sigma}_{ij,j}\, \E_i\otimes\E_j)\E_j\,.
\end{equation}
Finally, $\tilde{\sigma}_{ij,j}$ is computed  as:
\begin{equation}\label{eq:stress_div_global_comp}
\tilde{\sigma}_{ij,j}=\frac{\diff\tilde{\sigma}_{ij}}{\diff X_j}=\left(\E_i\otimes\E_j\right):\T{\nabla} \cdot \T{\sigma}=\left(\E_i\otimes\E_j\right):\T{\nabla} \cdot \CC : \strain +
\left(\E_i\otimes\E_j\right) : \CC \tensorm \T{\nabla}\strain\,.
\end{equation}
We refer readers to Appendix~\ref{sec:stress_comp} for further details
on the evaluation of these components. 

\subsection{Principle of virtual work}\label{subsec:weak_form}
The momentum balance equation can be imposed in a weak sense relying on the principle of virtual
work, which states that the sum of the system virtual internal, $\delta W_{int}$, and external, $\delta W_{ext}$, work is zero in an equilibrium state:
\begin{equation}\label{eq:PVW}
\delta W_{int} + \delta W_{ext} = 0\,.
\end{equation}
Therefore, the linear elasticity problem in variational form, reads
\begin{equation}\label{eq:weak_form}
\int_{\Omega}\boldsymbol{\sigma}:\delta\boldsymbol{\varepsilon}\diff\Omega-\int_{\Omega}\boldsymbol{b}\cdot\delta\boldsymbol{u}\diff\Omega-\int_{\Gamma_N}\boldsymbol{t}\cdot\delta\boldsymbol{u}\diff\Gamma=0\,,
\end{equation}
where $\delta\boldsymbol{u}$ and $\delta\boldsymbol{\varepsilon}$, respectively, are the virtual displacement field and strain defined in the global reference system $\{X_i\}$.
	
    \section{Stress recovery for curved laminated composite structures}\label{sec:recovery}
    In this section we describe the proposed post-processing strategy for curved structures, which is general and, provided the sufficient continuity, allows to recover the out-of-plane stresses independently of the designated numerical method to approximate the displacement field.

\begin{figure}
	\centering
	\includegraphics[width=0.5\textwidth]{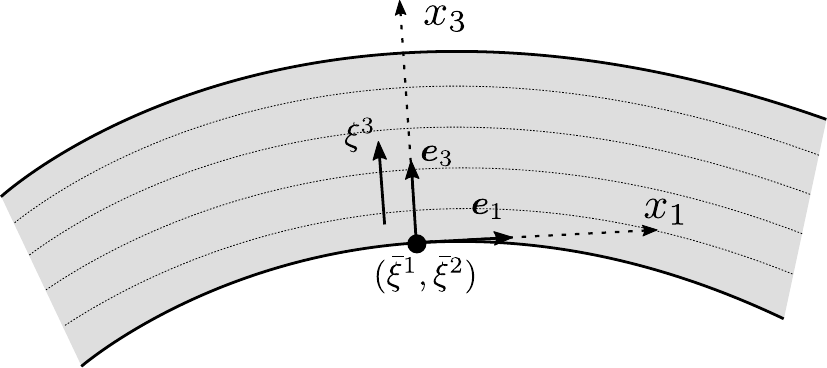}
	\caption{Local reference cartesian system at a point $(\bar\xi^1,\bar\xi^2)$
		for the stress recovery.}
	\label{fig:02}
\end{figure}

Thus, in an equilibrium state, stresses inside the material should satisfy at every point the equilibrium equation~\eqref{eq:equil}, which using the global reference system $\{X_i\}$ can be further detailed in a componentwise way as:
\begin{subequations}\label{eq:equilibrium_global}
\begin{align}
\tilde\sigma_{11,1} + \tilde\sigma_{12,2} + \tilde\sigma_{13,3} + \tilde b_1 &= 0\,,\label{eq:equilibrium_global_1}\\
\tilde\sigma_{12,1} + \tilde\sigma_{22,2} + \tilde\sigma_{23,3} + \tilde b_2&= 0\,,\label{eq:equilibrium_global_2}\\
\tilde\sigma_{13,1} + \tilde\sigma_{23,2} + \tilde\sigma_{33,3} +\tilde b_3 &= 0\,.\label{eq:equilibrium_global_3}
\end{align}
\end{subequations}
In an analogous way, the equilibrium equations can be written with respect to the local cartesian reference system $\{x_\alpha\}$, such that $x_1$, $x_2$, and $x_3$ are taken as the fibre, matrix, and normal directions, respectively, as:
\begin{subequations}\label{eq:equilibrium_local}
\begin{align}
\sigma_{11,1} + \sigma_{12,2} + \sigma_{13,3} + b_1&=  0\,,\label{eq:equilibrium_local_1}\\
\sigma_{12,1} + \sigma_{22,2} + \sigma_{23,3} + b_2&=  0\,,\label{eq:equilibrium_local_2}\\
\sigma_{13,1} + \sigma_{23,2} + \sigma_{33,3} + b_3&=  0\,,\label{eq:equilibrium_local_3}
\end{align}
\end{subequations}
where $b_\alpha = \tilde{b}_i\,C_{i\alpha}$ and $\sigma_{\alpha\beta,\mu}$ are the stress derivatives:
\begin{align}\label{eq:derstress_comp_local}
\sigma_{\alpha\beta,\mu} = \frac{\diff \sigma_{\alpha\beta}}{\diff x_\mu}\,.
\end{align}
The advantage of using the local system $\{x_\alpha\}$ (see Figure~\ref{fig:02}) lies in the fact that the stress components $\sigma_{11}$, $\sigma_{22}$, and $\sigma_{12}$ are now the in-plane components, that are well-known starting from the obtained global displacement field $\tilde u_i$ using Equation~\eqref{eq:stress_local_global_coefs}, while $\sigma_{13}$, $\sigma_{23}$, and $\sigma_{33}$ represent the out-of-plane stresses, which are instead not correctly captured as it has already been investigated in \cite{Dufour2018,Patton2019} for the 3D plate case (and as it will be shown in Figures~\ref{fig:LWI_IGAG_stress_pt1_1_l11_S20} and~\ref{fig:LWI_IGAC_stress_pt1_1_l11_S20} of Section~\ref{sec:numerical_tests}).
Our goal is to obtain a good description of the out-of-plane stress components along the composite thickness direction, extending the procedure detailed in \cite{Dufour2018} for the 3D plate case. 
It is important to remark that due to the fact that the stress divergence is computed with respect to $\{x_\alpha\}$, that is a fixed cartesian reference system, no additional terms appear in Equation~\eqref{eq:equilibrium_local} as it would be the case for Equation~\eqref{eq:equilibrium_global}.
\begin{remark}\label{remark:3}
Another possibility would be
to express $\stress$ in its curvilinear components (either covariant or contravariant)
and to develop the term $\T{\nabla} \cdot \T{\sigma} + \T{b} = \T{0}$ using one of those bases.
However, using curvilinear coordinates, new terms would arise in the componentwise equilibrium equations
and most importantly, these new terms would couple together
the in-plane and out-of-plane stress components, making impossible to apply the stress recovery procedure proposed in~\cite{Dufour2018} in a straightforward manner.
\end{remark}
Thus, integrating Equations~\eqref{eq:equilibrium_local_1} and~\eqref{eq:equilibrium_local_2} along the thickness direction $x_3$, the shear out-of-plane stress components
can be computed as:
\begin{subequations}\label{eq:recoveryShear}
\begin{align}
\label{eq:s13}
\sigma_{13}(x_3) &= -\int^{x_3}_{\underline{x}_3}\big(\sigma_{11,1}(\zeta) + \sigma_{12,2}(\zeta)+b_1(\zeta)\big)\diff \zeta + \sigma_{13}(\underline{x}_3)\,,\\
\label{eq:s23}
\sigma_{23}(x_3) &= -\int^{x_3}_{\underline{x}_3}\big(\sigma_{12,1}(\zeta) + \sigma_{22,2}(\zeta)+b_2(\zeta)\big)\diff \zeta + \sigma_{23}(\underline{x}_3)\,,
\end{align}
\end{subequations}
where $\underline{x}_3$ indicates the value of $x_{3}$ at the bottom of the shell.

To recover also the out-of-plane normal stress profile $\sigma_{33}$, we focus on each $k$-th layer and substitute the appropriate derivatives of the out-of-plane shear stresses \eqref{eq:recoveryShear} into Equation~\eqref{eq:equilibrium_local_3}:
\begin{equation}\label{eq:s33,3}
\begin{aligned}
	\sigma^{(k)}_{33,3}(x^{(k)}_3) & = \int^{x^{(k)}_3}_{\underline{x}^{(k)}_3}\big(\sigma_{11,11}(\zeta) + \sigma_{22,22}(\zeta) + 2\sigma_{12,12}(\zeta)- b_{1,1}(\zeta)-b_{2,2}(\zeta)\big)\diff \zeta\\
	&-\big(\sigma_{13,1}(\underline{x}^{(k)}_3) + \sigma_{23,2}(\underline{x}^{(k)}_3)\big)-b_3(x^{(k)}_3)\,,
\end{aligned}
\end{equation}
where $\underline{x}^{(k)}_3\leq x^{(k)}_{3}\leq \bar{x}^{(k)}_3$, being $\underline{x}^{(k)}_3$ and $\bar{x}^{(k)}_3$, respectively, the values of the out-of-plane coordinate at the bottom and the top of the $k$-th layer, while
$\sigma_{\alpha\beta,\mu\nu}$ is the stress second derivative:
\begin{align}\label{eq:der2stress_comp_local}
\sigma_{\alpha\beta,\mu\nu} = \frac{\diff^2 \sigma_{\alpha\beta}}{\diff x_\mu \diff x_\nu}\,.
\end{align}
Finally, we further integrate Equation~\eqref{eq:s33,3} along the thickness obtaining
\begin{equation}\label{eq:s33}
 \sigma_{33}(x_3) = \int^{x_3}_{\underline{x}_3}\sigma_{33,3}(\xi)\diff\xi 
+ \sigma_{33}(\underline{x}_3)\,,
\end{equation}
where $\sigma_{33,3}(\xi)$ is equal to $\sigma^{(k)}_{33,3}(x^{(k)}_3)$ once it has been evaluated for every $k$-th layer.

For the sake of readability, the stress derivatives $\sigma_{\alpha\beta,\mu}$ and $\sigma_{\alpha\beta,\mu\nu}$ involved in Equations~\eqref{eq:recoveryShear} and~\eqref{eq:s33,3} are detailed in Appendix~\ref{sec:stress_comp}.

\begin{remark}\label{remark:4}
It should be noted that integrals~\eqref{eq:recoveryShear},~\eqref{eq:s33,3}, and~\eqref{eq:s33} are computed numerically using a composite trapezoidal quadrature rule. Also, new stress boundary conditions (namely, $\sigma_{13,1}(\underline{x}^{(k)}_3)$ and $\sigma_{23,2}(\underline{x}^{(k)}_3)$) arise in Equation~\eqref{eq:s33,3} for each $k$-th layer. 
\end{remark}

With reference to Equations~\eqref{eq:recoveryShear},~\eqref{eq:s33,3}, and~\eqref{eq:s33}, it is clear that to apply the proposed post-processing technique a highly regular displacement solution is needed. More specifically, the required in-plane derivatives $\sigma_{\alpha\beta,\mu}$ and $\sigma_{\alpha\beta,\mu\nu}$, detailed respectively in the Appendix in Equation~\eqref{eq:derstress_local_comp} and~\eqref{eq:der2stress_local_comp}, need to be computed from a  C$^2$-continuous in-plane displacement solution, which can be achieved by means of, \eg isogeometric analysis.
	
	\section{IGA strategies for 3D laminated curved geometries: Basics and application to structures made of multiple orthotropic layers}\label{sec:numerical_strategies}
	In Section \ref{sec:recovery} we described our stress recovery technique which can be applied regardless the numerical approximation method, as long as able to provide the required continuity. Due to its high-order continuity properties, IGA proves to be a natural choice for the proposed approach, while other techniques, such as the Finite Element method, would require a hybrid a-posteriori interpolation via, \eg splines or radial basis functions.

Thus, we hereby introduce the notions of multivariate NURBS, and detail the proposed displacement-based IGA Galerkin and Collocation strategies to analyze laminated composite curved geometries. 

\subsection{Multivariate NURBS}\label{subsec:IGAmultivarNURBS}
In the following, we introduce the basic definitions and notations about multivariate NURBS. For further details, readers may refer to \cite{Cottrell2007,Hughes2005,Piegl1997}, and references therein. 
Having defined $d_{s}$ as the dimension of the considered physical space, NURBS geometries in $\mathbb{R}^{d_{s}}$ are constructed starting from a projective transformation of their B-Spline counterparts in $\mathbb{R}^{d_{s}+1}$.
Therefore, to introduce the notion of multivariate B-Splines, which are obtained through the tensor product of univariate B-Splines, a few basic quantities need to be defined.
We denote with $d_{p}$ the dimension of the parameter space and $d_{p}$
univariate knot vectors, \ie non-decreasing set of coordinates in the $d$-th parameter
space, can be introduced as
\begin{alignat}{2}\label{eq:knotvectors}
&\Xi^{d} = \{\xi_{1}^{d},...,\xi_{m_{d}+p_{d}+1}^{d}\} \quad &\quad d = 1, ..., d_{p}\,,
\end{alignat}
where $p_{d}$ represents the polynomial degree in the parametric direction $d$ (such that in the case of $d=1,2,3$, for the sake of clarity, we name the in-plane degree of approximation $p_1=p$, $p_2=q$, and the out-of-plane one $p_3=r$),
while $m_{d}$ is the associated number of basis functions. A univariate B-Spline basis function $N^{d}_{i_{d},p_{d}}(\xi^{d})$ associated to each parametric coordinate $\xi^{d}$ can be then constructed for each $i_d$ position in the tensor
product structure, using the Cox-de Boor recursion formula, starting from $p_{d}=0$, as
\begin{alignat}{2}\label{eq:univarBSplinesConst}
&N^{d}_{i_{d},0}(\xi^{d})=\begin{cases}1\quad &\quad\xi_{i_{d}}^{d}\le\xi^{d}<\xi_{i_{d}+1}^{d}\\
0\quad&\quad \text{otherwise}
\end{cases}\,,
\end{alignat}
while the basis functions for $p_{d}>0$ are obtained from
\begin{equation}\label{eq:univarBSplines}
N^{d}_{i_{d},p_{d}}(\xi^{d})=\cfrac{\xi^{d}-\xi^{d}_{i_{d}}}{\xi^{d}_{i_{d}+p_{d}}-\xi^{d}_{i_{d}}}N^{d}_{i_{d},p_{d}-1}(\xi^{d})+\cfrac{\xi^{d}_{i_{d}+p_d+1}-\xi^{d}}{\xi^{d}_{i_{d}+p_d+1}-\xi^{d}_{i_{d}+1}}N^{d}_{i_{{d}+1},p_{d}-1}(\xi^{d}),
\end{equation}
where the convention $0/0=0$ is assumed.
Having defined the univariate basis functions $N^{d}_{i_{d},p_{d}}(\xi^{d})$, the multivariate basis functions $B_{\textbf{i},\textbf{p}}(\boldsymbol{\xi})$ are obtained as:
\begin{equation}\label{eq:multivarBSplines}
B_{\textbf{i},\textbf{p}}(\boldsymbol{\xi})=\prod\limits_{d=1}^{d_{p}}N_{i_{d},p_{d}}(\xi^{d})\,,
\end{equation}
where $\textbf{i} = \{i_{1}, ..., i_{d_{p}}\}$ plays the role of a multi-index which describes the considered position in the tensor product structure, $\textbf{p} = \{{p_{1},..., p_{d}}\}$ indicates the polynomial degrees, and $\boldsymbol{\xi} = \{\xi^{1},...,\xi^{d_{p}}\}$ represents the vector of the parametric coordinates in each parametric direction $d$. B-Spline multidimensional geometries are built from a linear combination of multivariate B-Spline basis functions as follows
\begin{equation}\label{eq:multivarBSplinesGeom}
\textbf{S}(\boldsymbol{\xi})=\sum\limits_{\textbf{i}}B_{\textbf{i},\textbf{p}}(\boldsymbol{\xi})\textbf{P}_\textbf{i}\,,
\end{equation}
where the coefficients $\textbf{P}_\textbf{i}\in\mathbb{R}^{d_{s}}$ are the so-called
control points and the summation is extended to all combinations of the multi-index \textbf{i}.
Finally we can introduce multivariate NURBS basis functions as
\begin{equation}\label{eq:multivarNURBS}
R_{\textbf{i},\textbf{p}}(\boldsymbol{\xi})=\cfrac{B_{\textbf{i},\textbf{p}}(\boldsymbol{\xi})w_{\textbf{i}}}{\sum_{\textbf{j}}B_{\textbf{j},\textbf{p}}(\boldsymbol{\xi})w_{\textbf{j}}}\,,
\end{equation}
where $w_{\textbf{i}}$ represent the collection of NURBS weights associated to each control point according to the multi-index \textbf{i}.
NURBS multidimensional geometries are then built combining multivariate NURBS basis functions and control points as follows  
\begin{equation}
\textbf{S}(\boldsymbol{\xi})=\sum\limits_{\textbf{i}}R_{\textbf{i},\textbf{p}}(\boldsymbol{\xi})\textbf{P}_\textbf{i}\,.\label{eq:multivarNURBSgeom}
\end{equation}

\subsection{Numerical strategies}\label{subsec:numerical_strategies}
The considered modeling strategy can be regarded as a three-dimensional equivalent-single-layer (ESL) approach, which models the 3D laminate employing only one element through the thickness, which strongly reduces the
number of degrees of freedom with respect to layerwise methods.
Besides the intrinsic nature of IGA-Galerkin and Collocation methods in terms of how the linear balance momentum equation is approximate (either at the strong or at the weak form level), the two approaches differ also in the way material properties are considered. 
For an IGA Galerkin method,
constitutive features can be taken into account using through-the-thickness either a calibrated layerwise integration rule or a homogenized approach, with the latter representing the most direct ESL strategy for Collocation to the authors' knowledge.

\subsubsection{IGA-Galerkin method strategy}\label{subsubsec:galerkin}
Addressing the element point of view and adopting a standard Gauss quadrature, the global displacement $\boldsymbol{u}$ and virtual displacement field $\delta\boldsymbol{u}$ are approximated as a linear combination of multivariate shape functions $R^{(e)}_{\textbf{i},\textbf{p}}(\boldsymbol{\bar{\xi}})$ (defined in Equation~\eqref{eq:multivarNURBS}) and control variables $\hat\uu^{(e)}_{\textbf{i}}$ as 
\begin{subequations}\label{eq:displapproxGal}
	\begin{align}
	&\bm{u}^{(e)}\simeq\bm{u}^{(e)}_h(\boldsymbol{\bar{\xi}})=\RR^{(e)}_{\textbf{i},\textbf{p}}(\boldsymbol{\bar{\xi}})\,\hat\uu^{(e)}_{\textbf{i}}\,,\label{eq:displapproxGalu}\\
	&\delta\bm{u}^{(e)}\simeq\delta\bm{u}^{(e)}_{h}(\boldsymbol{\bar{\xi}})=\RR^{(e)}_{\textbf{i},\textbf{p}}(\boldsymbol{\bar{\xi}})\,\delta\hat\uu^{(e)}_{\textbf{i}}\,,\label{eq:displapproxGaldeltaw}
	\end{align}
\end{subequations}
where the superscript $(e)$ is the element index, while with $h$ we denote any approximated field.
The multivariate shape function element matrix $\RR^{(e)}_{\textbf{i},\textbf{p}}(\boldsymbol{\bar{\xi}})$ is instead defined as
\begin{equation}\label{eq:sf_matrix_el}
\RR^{(e)}_{\textbf{i},\textbf{p}}(\boldsymbol{\bar{\xi}})=
\begin{bmatrix}	
R^{(e)}_{\textbf{i},\textbf{p}}(\boldsymbol{\bar{\xi}})&0&0\\
0&R^{(e)}_{\textbf{i},\textbf{p}}(\boldsymbol{\bar{\xi}})&0\\
0&0&R^{(e)}_{\textbf{i},\textbf{p}}(\boldsymbol{\bar{\xi}})
\end{bmatrix}\,,
\end{equation}
being $\boldsymbol{\bar{\xi}}$ the matrix of the quadrature point positions which are of relevance for the considered element $e$.

\begin{figure}[!htbp]
	\centering
	\scalebox{1.}{{\ifrecompiletikz\tikzsetnextfilename{fig_03}\tikzexternalenable\input{images/fig_03}\tikzexternaldisable\else\includegraphics{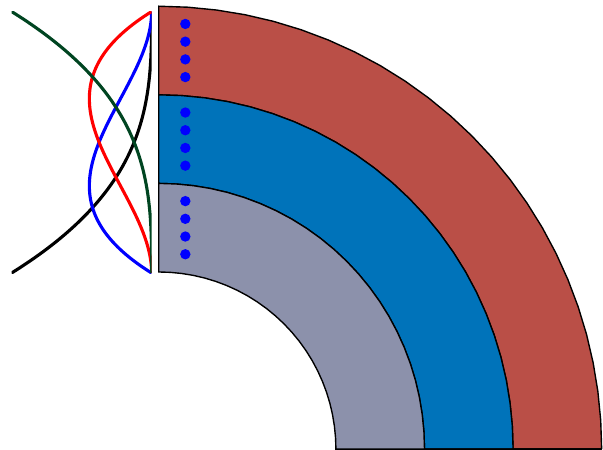}\fi}}
	\caption{Single-element approach for IGA Galerkin method with special through-the-thickness integration rules (\ie $r + 1$ Gauss points per layer). Example of shape functions for an out-of-plane degree of approximation $r = 3$. The blue bullets represent the position of the quadrature points along the thickness.}
	\label{fig:Gal_LSE}
\end{figure}

Thus, we discretize Equation~\eqref{eq:weak_form} obtaining the approximated element energy variation as 
\begin{equation}\label{eq:weak_form_global_discr}
\delta W_{int}+\delta W_{ext}=\sum_{e=1}^{nel}\int_{\Omega^{(e)}}(\CC^{(k)}:\T{\nabla}^{s}\bm{u}^{(e)}_h(\boldsymbol{\bar{\xi}})):\T{\nabla}^s\delta\bm{u}^{(e)}_h(\boldsymbol{\bar{\xi}})\,\rmd\Omega^{(e)}\,.
\end{equation}
Laminated composites often exhibit material properties which may vary layer by layer even when they are pointwisely referred to the principal material coordinates.
Using a special plywise integration rule (namely, $r+1$ Gauss point per layer as displayed in Figure~\ref{fig:Gal_LSE}), the Galerkin method provides a natural approach to account for this through-the-thickness dependency.
In fact, we remark that $\CC^{(k)}$ (\ie $\CC^{(k)}=\tilde\CC^{(k)}_{ijkl}\,\E_i\otimes\E_j\otimes\E_k\otimes\E_l$, with $\tilde\CC^{(k)}_{ijkl}$ defined in Equation~\eqref{eq:CC_change_basis}) represents the global material property tensor for the $k$-th layer, which can be traced during the elementwise assembly relying on the out-of-plane Gauss point number, allowing to significantly improve the overall post-processing quality as it will be clear in Section~\ref{sec:numerical_tests}.
\begin{remark}\label{remark:1}
	Using a layerwise integration rule allows to correctly capture the behavior of composites for any stacking sequence (\eg both an even and odd number of variously oriented layers can be accounted for), or consider more general constitutive models, such as plasticity. Furthermore, we highlight that, based on our numerical experiments, considering $r + 1$ quadrature points per layer leads to basically the same accuracy as using $r-1$, in accordance with the solid plate case (see \cite{Dufour2018}).
	
	An alternative to the proposed integration method would be to split the assembly of stiffness matrices into their in-plane and out-of-plane contributions as in \cite{Antolin2020a}, reducing the assembly cost significantly. 
\end{remark}

\subsection{IGA-Collocation strategy}\label{subsubsec:collocation}
\begin{figure}[!htbp]
	\centering
    \scalebox{1.0}{{\ifrecompiletikz\tikzsetnextfilename{fig_04}\tikzexternalenable\input{images/fig_04}\tikzexternaldisable\else\includegraphics{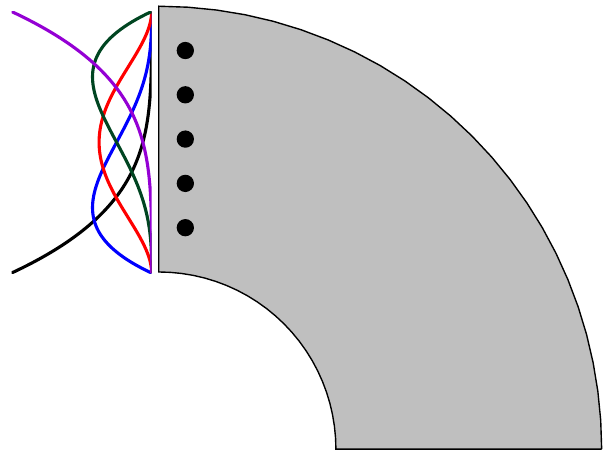}\fi}}
	\caption{Homogenized single-element approach for isogeometric collocation ($r + 1$ evaluation points independently on the number of layers). Example of shape functions for an out-of-plane degree of approximation $r = 4$. The black bullets represent the position of the quadrature points along the thickness.}
	\label{fig:Coll_HSE}
\end{figure}

In this section we describe our IGA 3D collocation strategy to model laminated composite curved structures.
Following \cite{Patton2019}, we propose a homogenized single element approach, which takes into account layerwise variations of orthotropic material properties, homogenizing the constitutive behavior to form an equivalent single-layer laminate as in Figure~\ref{fig:Coll_HSE}. While this approach represents an obvious choice in the context of collocation, it can be regarded instead as a less accurate but cheaper alternative for the IGA-Galerkin method proposed in Section~\ref{subsubsec:galerkin}. In fact, to compute the displacement solution needed for the post-processing, the homogenized approach requires $r+1$ integration points regardless the number of layers, while, using the previously introduced special integration rule, $r+1$ Gauss points are employed for each ply.

Under these premises, for collocation, we homogenize the material properties according to \cite{Sun1988} by means of the following relations
\begin{subequations}\label{eq:aveC}
	\begin{alignat}{2}
	&\CCo_{\zeta\eta}=\sum_{k=1}^{N}\overline{t}_{k}\CC_{\zeta\eta}^{(k)}+\sum_{k=2}^{N}(\CC_{\zeta 3}^{(k)}-\CCo_{\zeta 3})\overline{t}_{k}\frac{(\CC_{\eta 3}^{(1)}-\CC_{\eta 3}^{(k)})}{\CC_{33}^{(k)}}&\quad\quad\begin{split}
	\begin{aligned}&zeta=1,2\\&\eta=1,2,3\end{aligned}
	\end{split}\,,\label{eq:aveC1}\\
	&\CCo_{33}=\frac{1}{\bigg(\sum_{k=1}^{N}\cfrac{\overline{t}_{k}}{\CC_{33}^{(k)}}\bigg)}\,,&\quad\quad&\label{eq:aveC2}\\
	&\CCo_{\theta\theta}=\frac{\bigg(\sum_{k=1}^{N}\cfrac{\overline{t}_{k}\CC_{\theta\theta}^{(k)}}{\Delta_{k}}\bigg)}{\Delta}
&\quad\quad& \theta=4,5\,,\label{eq:aveC3}\\
	&\CCo_{66}=\sum_{k=1}^{N}\overline{t}_{k}\CC_{66}^{(k)}\,,&\quad\quad&\label{eq:aveC9}
	\end{alignat}
\end{subequations}
with $\Delta=\prod_{\theta\theta=4}^5\bigg(\sum_{k=1}^{N}\cfrac{\overline{t}_{k}\CC_{\theta\theta}^{(k)}}{\Delta_{k}}\bigg)$ and $\Delta_{k}=\prod_{\theta\theta=4}^5\CC_{\theta\theta}^{k}\,$. 

In Equation~\eqref{eq:aveC}, $\CC_{\zeta\eta}^{(k)}$ represents the $\zeta\eta$-th local component of the fourth order elasticity tensor in Voigt' notation~\eqref{eq:CC_coeffs} for the $k$-th layer and $\overline{t}_{k}=t_{k}/h$ stands for the volume fraction of the $k$-th lamina, $h$ being the total thickness and $t_{k}$ the $k$-th thickness.
In order to be used in the global framework (see Equation~\eqref{eq:stress_div_global_comp}), the homogenized material tensor $\CCo=\CCo_{\alpha\beta\gamma\delta}\,\a_\alpha\otimes\a_\beta\otimes\a_\gamma\otimes\a_\delta$, mapped from Voigt's to its full representation in indicial notation, is transformed according to the basis change operator $D_{i\alpha}$, analogously to Equation~\eqref{eq:CC_change_basis}, as 
\begin{equation}\label{eq:CC_rot_homo}
\tilde{\overline{\CC}}_{ijkl}=D_{i\alpha}\,D_{j\beta}\,D_{k\gamma}\,D_{l\delta}\,\CCo_{\alpha\beta\gamma\delta}\,.
\end{equation}
\begin{remark}\label{remark:2}
	Considering a local reference system allows to define whether the ply stacking sequence is symmetric or not with respect to the solid geometrical mid-surface.  
	Adopting this pointwise perspective, we remark that the homogenized single-element approach is immediately effective only for symmetric ply distributions as, for non-symmetric ones, the laminate geometric and material mid-surface do not coincide. In any case, symmetric stacking sequences typically cover the most common cases in practice and, in the need for laminates made of non-symmetric layer distributions, this technique finds still application when the stacking sequence can be split into two symmetric piles, using one element per homogenized stack with a $C^0$ interface \cite{Patton2019}.
\end{remark}
Having defined how material properties are tuned for the homogenized single-element approach, we procede to detail our collocation strategy. Collocation methods (see \cite{Auricchio2010} for further details) directly discretize in strong form the differential equations governing the problem evaluated at collocation points. Therefore, a delicate issue lies in how these collocation points are determined. The simplest and most widespread approach in the engineering literature (see \cite{deBoor1973,Demko1985} for alternative
choices) is
to collocate at the images of Greville abscissae, \ie points obtained from the knot vector components, $\xi^d_{i}$, as
\begin{alignat}{2}\label{eq:greville}
&\bar{\xi}^d_{i_d}=\frac{\xi^d_{i_d+1}+\xi^d_{i_d+2}+...+\xi^d_{i_d+p_d}}{p_d}\quad&&\quad i_d = 1,...,m_d\,.
\end{alignat}

Then, having defined $\boldsymbol{\tau}$ as the matrix of collocation points, the global displacement field $\boldsymbol{u}$ is approximated as a linear combination of NURBS multivariate shape functions $R^{(e)}_{\textbf{i},\textbf{p}}(\boldsymbol{\tau})$ and control variables $\hat{\uu}_{\textbf{i}}$ as
\begin{equation}\label{eq:displapproxCol}
\boldsymbol{u}\simeq \boldsymbol{u}_{h}(\boldsymbol{\tau})=\tilde\RR_{\textbf{i},\textbf{p}}(\boldsymbol{\tau})\,\hat\uu_{\textbf{i}}\,,
\end{equation} 
where the $\boldsymbol{\tau}$ matrix has been defined such that for each $i_d$-th point and $d$-th parameter space $\tau^{d}_{i_d}=\cfrac{\sum_{l=1}^{p_d}\xi_{i_d+l}}{p_d}$ with $i_d = 1,...,m_d$.
The multivariate shape functions matrix $\tilde\RR_{\textbf{i},\textbf{p}}(\boldsymbol{\tau})$ is instead characterized as follows
\begin{equation}\label{eq:sf_matrix}
\RR_{\textbf{i},\textbf{p}}(\boldsymbol{\tau})=
\begin{bmatrix}	
R_{\textbf{i},\textbf{p}}(\boldsymbol{\tau})&0&0\\
0&R_{\textbf{i},\textbf{p}}(\boldsymbol{\tau})&0\\
0&0&R_{\textbf{i},\textbf{p}}(\boldsymbol{\tau})
\end{bmatrix}\,.
\end{equation}
Finally, the approximation of system~\eqref{eq:strong_form} reads as
 \begin{subequations}\label{eq:strong_form_approx}
	\begin{alignat}{3}
	&\T{\nabla} \cdot \CCo : \T{\nabla}^{s}\bm{u}_h(\tau^{d}_{i}) + \CCo \tensorm \boldsymbol{\nabla}(\T{\nabla}^{s}\bm{u}_h(\tau^{d}_{i})) + \T{b}_h(\tau^{d}_{i}) = \T{0}                                   \quad\quad &    \quad & \quad\quad\forall\tau^{d}_{i}\in\Omega\,,
	\label{eq:equil_approx}\\
	& \CCo : \T{\nabla}^{s}\bm{u}_h(\tau^{d}_{i}) \cdot \boldsymbol{n}(\tau^{d}_{i})=\boldsymbol{t}_h(\tau^{d}_{i})   \quad\quad &    \quad & \quad\quad\forall\tau^{d}_{i}\in\Gamma_N\,,\label{eq:nbc_approx}\\
	& \boldsymbol{u}_h(\tau^{d}_{i}) = \bar{\boldsymbol{u}}_h(\tau^{d}_{i})                                                   \quad\quad &    \quad & \quad\quad\forall\tau^{d}_{i}\in\Gamma_D\,,\label{eq:dbc_approx}
	\end{alignat}
\end{subequations}
where all quantities are evaluated in the global reference system as detailed in Section~\ref{sec:governing_eqq}.
     
	\section{Numerical Tests}\label{sec:numerical_tests}
	In this section, we propose several benchmarks to showcase the effectiveness of the presented approach. To this extent, we consider a composite solid hollow cylinder under bending and we validate the obtained results against an overkilled $C^0$ layerwise solution, addressing the main differences with respect to the proposed through-the-thickness integration IGA Galerkin and homogenized IGA collocation approaches, as well as the method sensitivity to parameters of interest (\ie number of layers and thickness-to-mean radius ratio).

\subsection{Composite solid cylinder under bending}\label{subsec:cylinder_description}
\begin{figure}[htbp!]
	\centering
	\scalebox{0.7}{\ifrecompiletikz\tikzsetnextfilename{fig_05}\tikzexternalenable\input{images/fig_05}\tikzexternaldisable\else\includegraphics{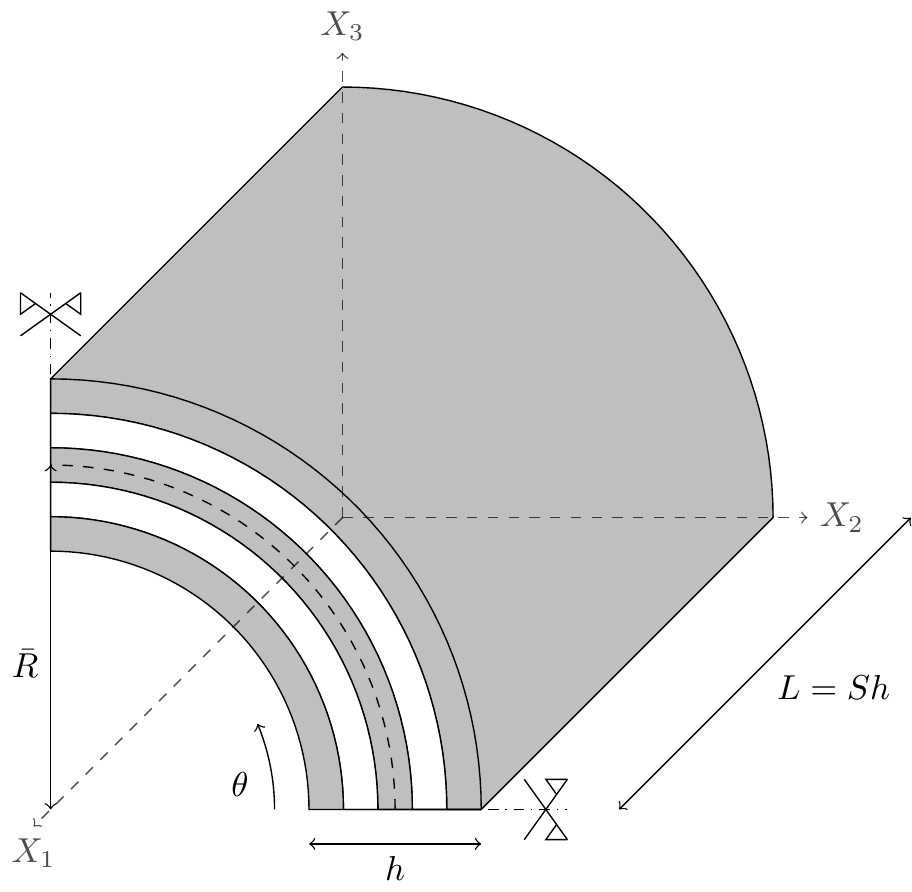}\fi}
	\caption{Quarter of composite cylindrical shell: Problem geometry.}
	\label{fig:quarter_cylinder}
\end{figure}

A solid hollow cross-ply cylindrical shell of total thickness $h$ and made of $N$ orthotropic layers is considered as in \cite{Varadan1991}. The structure, as detailed in Figure~\ref{fig:quarter_cylinder}, is simply supported at both ends and subjected to a transverse sinusoidal loading, $q(X_1,X_2,X_3)$ on the inner surface, while the top surface is traction-free. 
\newline
The thickness of every single layer is set to 1 mm and the edge length $L$ is chosen to be equal to the mean radius $\bar{R}$, such that $\bar{R}$ is $S$ times larger than the total thickness of the laminate, $S=\bar{R}/h=L/h$. Thus, $S$ represents the inverse of the slender parameter or thickness-to-mean radius ratio.  
\begin{table}[!htbp]
	\caption{Numerical tests material properties for 0\textdegree-oriented layers.} 
	\begin{center}
		\small{
			\begin{tabular}{?c c c c c c c c c?}
				\thickhline 
				$E_{1}$ & $E_{2}$& $E_{3}$ & $G_{23}$ & $G_{13}$ & $G_{12}$ & $\nu_{23}$ & $\nu_{13}$ & $\nu_{12}$\Tstrut\Bstrut\\\hline 
				[MPa] & [MPa] & [MPa] & [MPa] & [MPa] & [MPa]& [-] & [-] & [-]\Tstrut\Bstrut\\\thickhline 
				25 & 1 & 1 & 0.2 & 0.5 & 0.5 & 0.25 & 0.25& 0.25\Tstrut\Bstrut\\\thickhline  
		\end{tabular}}
	\end{center}\label{tab:matproperties}
\end{table}

Layer material parameters taken into account for all the proposed numerical tests are summarized in Table~\ref{tab:matproperties} for 0\textdegree-oriented plies, while the considered loading pressure is equal to
\begin{equation}\label{eq:load}
q(X_1,X_2,X_3)=\sigma_0\cos(4\theta)\sin(\dfrac{\pi X_1}{Sh})\,,
\end{equation}
where $\sigma_0 =$ -1 MPa and $\theta=\theta(X_2,X_3)\in[0,\bar{\theta}]=[0,\pi/2]$ is oriented anticlockwise as depicted in Figure~\ref{fig:quarter_cylinder}.
The inner and the outer surface radii are, respectively, $r_i=\bar{R}-\dfrac{h}{2}$ and $r_o=\bar{R}+\dfrac{h}{2}$.

The simply supported edge conditions are taken as
\begin{equation}
\begin{aligned}
\tilde\sigma_{11}=0\;\text{and}\;\tilde u_2=\tilde u_3=0\;\text{at}\;X_1=0\;\text{and}\;X_1=L\,,\label{eq:DBCs}
\end{aligned}
\end{equation}
while Neumann boundary conditions on the tube inner and outer surfaces are
\begin{subequations}\label{eq:NBCs}
	\begin{alignat}{2}
	&\tilde\sigma_{13}=\tilde\sigma_{23}=0 \;\text{and}\;\tilde\sigma_{33}=q(X_1,X_2,X_3)\quad&&    \quad\text{at}\;\sqrt{X_2^2+X_3^2}=r_i\,,\label{eq:NBCs_1}\\
	&\tilde\sigma_{13}=\tilde\sigma_{23}=\tilde\sigma_{33}=0\quad&&    \quad\text{at}\;\sqrt{X_2^2+X_3^2}=r_o\,.\label{eq:NBCs_2}
	\end{alignat}
\end{subequations}
Taking advantage of the problem's symmetry, we model only a quarter of  cylinder, adding the following further constraints
\begin{subequations}\label{eq:sym_BCs}
	\begin{align}
	&\tilde\sigma_{23}=\tilde\sigma_{12}=0\;\text{and}\;\tilde u_2=0\;\text{at}\;X_2=0\;,\label{eq:sym_BCs_1}\\
	&\tilde\sigma_{23}=\tilde\sigma_{13}=0\;\text{and}\;\tilde u_3=0\;\text{at}\;X_3=0\;.\label{eq:sym_BCs_2}
	\end{align}
\end{subequations}

All results in the local reference system are then normalized as
\begin{alignat}{2}\label{eq:normalizedresults}
\overline{\sigma}_{ij}=\dfrac{\sigma_{ij}}{\lvert\sigma_0\rvert}\quad&\quad i,j=1,2,3\,.
\end{alignat}

\subsection{Single-element approach results: The post-processing effect}\label{subsec:layerwise}
In this section, we present and comment several numerical results which consider a quarter of composite cylindrical shell with $S=20$ and a cross-ply distribution of 11 layers, namely a 0\textdegree/90\textdegree~stacking sequence. Layerwise methods allow to capture more accurately the mechanical state inside the laminate. Therefore, the stress profiles, either obtained via the appropriate constitutive law or recovered through the equilibrium imposition, are validated against an overkill IGA layerwise solution, computed according to \cite{Dufour2018}, which uses a degree of approximation $p=q=6$, $r=4$ and a number of control points equal to 36x36x5x$n_l$, where $n_l$ is the number of layers. More specifically, this layerwise reference solution is computed modeling each material layer of the laminate by one patch through the thickness and keeping $C^0$ continuity at each ply interface. This technique, despite being computationally demanding since it requires a number of unknowns directly proportional to the number of layers, allows to obtain a continuous out-of-plane stress profile through the thickness.

\subsubsection{IGA-Galerkin method with an ad-hoc through-the-thickness integration rule}
The presented numerical simulations are post-processed starting from a single-element displacement-based IGA-Galerkin approximation with an ad-hoc integration rule which uses $r+1$ Gauss points per layer and are obtained employing 22x22x4 control points, an in-plane degree of approximation $p=q=4$, and an out-of plane one $r=3$.

As a first example, in Figure~\ref{fig:LWI_IGAG_stress_pt1_1_l11_S20}, we present in the left column solution profiles for a sampling point located at ($X_1=L/3$, $\theta=\bar{\theta}/3$), proving the ability of the proposed approach to grant accurate in-plane results.
Also, in Figure~\ref{fig:LWI_IGAG_stress_pt1_1_l11_S20} the right column compares the out-of-plane stress state for the same sampling point with and without applying the presented post-processing step (see Section~\ref{sec:recovery}), showing a remarkable improvement in the solution profiles which, after recovery, are continuous through-the-thickness as required by equilibrium. In more detail, the use of a single element through the thickness with $C^\infty$ shape functions leads to continuous through-the-thickness displacements and hence strains. Then, if the latters are multiplied by layerwise discontinuous material properties, as in the case of cross-ply laminates, this implies a numerical approximation of the  stresses that is discontinuous through the thickness, differently from what prescribed by equilibrium for out-of-plane stress components.
%
\begin{figure}[!htbp]
	\centering
	\subfigure[Normalized $\sigma_{11}$\label{subfig-1:test11}]{\ifrecompiletikz\tikzsetnextfilename{fig_06_a}\tikzexternalenable\input{images/fig_06_a}\tikzexternaldisable\else\includegraphics{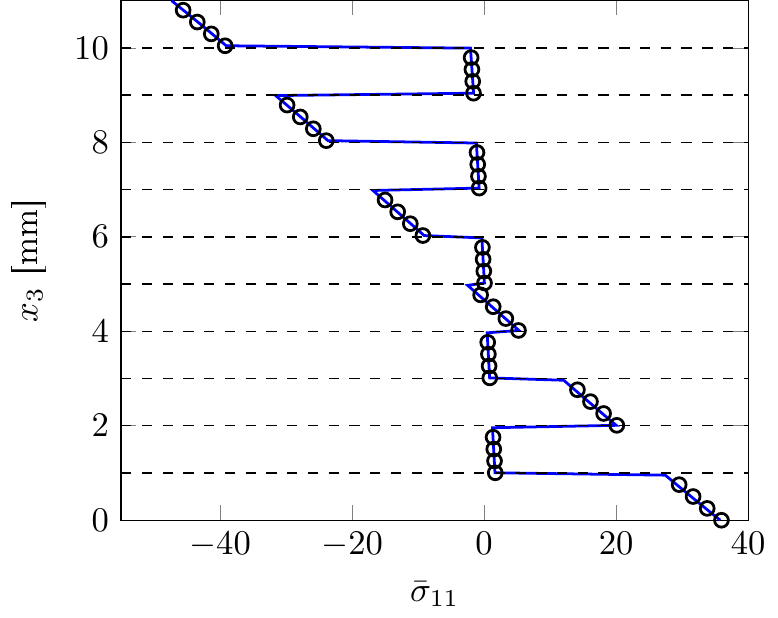}\fi}\hfill
	\subfigure[Normalized $\sigma_{13}$\label{subfig-2:test11}]{\ifrecompiletikz\tikzsetnextfilename{fig_06_d}\tikzexternalenable\input{images/fig_06_d}\tikzexternaldisable\else\includegraphics{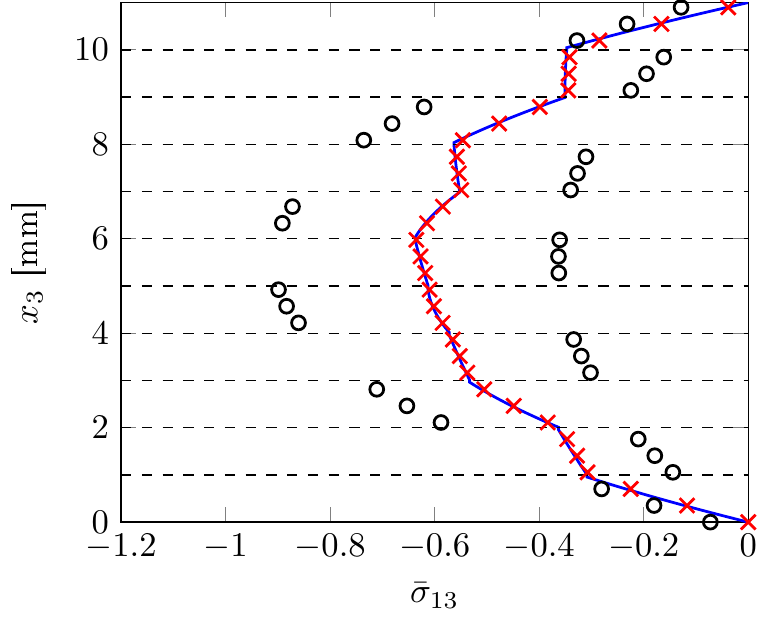}\fi}
	\caption{Through-the-thickness stress profiles evaluated at ($X_1=L/3$, $\theta=\bar{\theta}/3$) for IGA-Galerkin (degree of approximation $p=q=4$, $r=3$, and 22x22x4 control points). Case: Tube with mean radius-to-thickness ratio $S=20$, 11 layers, and $L=\bar{R}$ (\blueline~overkill IGA layerwise solution,~\blackcircles~single-element approach solution without post-processing, \redcrosses~post-processed solution).}
	\end{figure}
	\begin{figure}[!htbp]
	\centering
	\ContinuedFloat
	\captionsetup{list=off,format=cont}
	\subfigure[Normalized $\sigma_{22}$\label{subfig-3:test11}]{\ifrecompiletikz\tikzsetnextfilename{fig_06_b}\tikzexternalenable\input{images/fig_06_b}\tikzexternaldisable\else\includegraphics{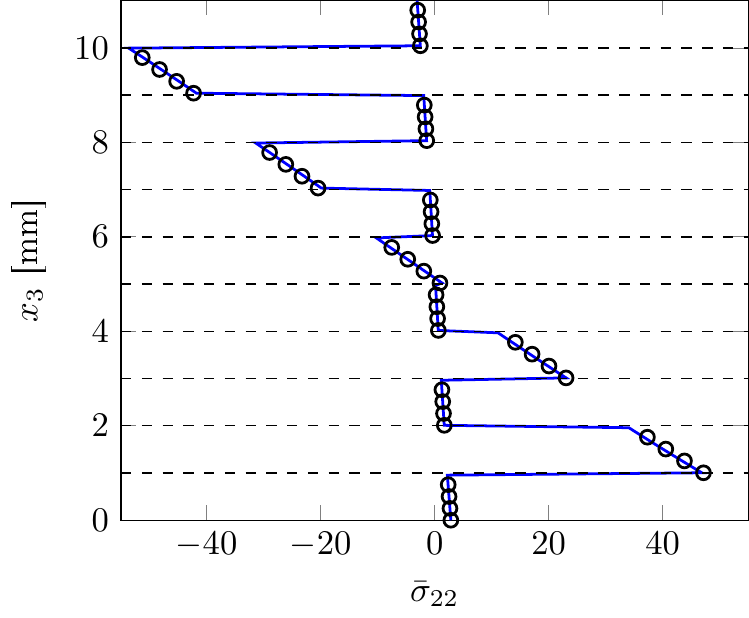}\fi}\hfill
	\subfigure[Normalized $\sigma_{23}$\label{subfig-4:test11}]{\ifrecompiletikz\tikzsetnextfilename{fig_06_e}\tikzexternalenable\input{images/fig_06_e}\tikzexternaldisable\else\includegraphics{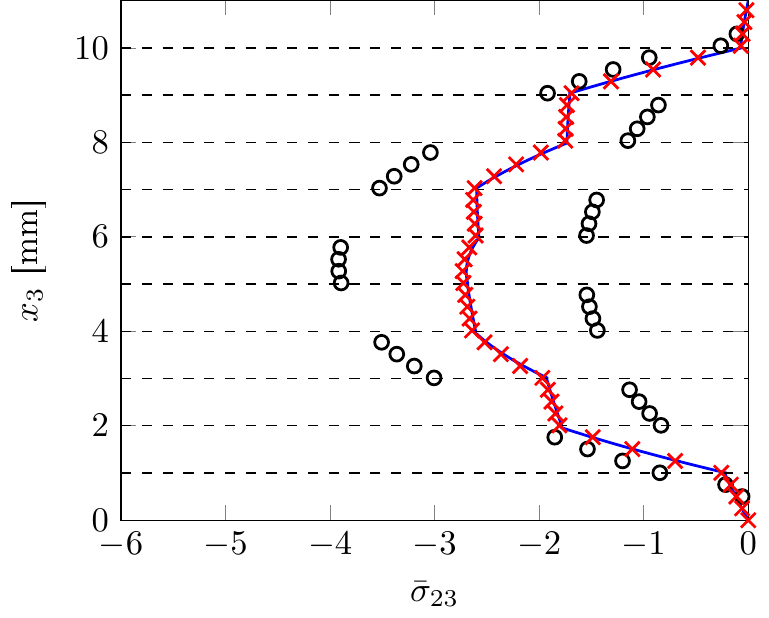}\fi}\\
	\subfigure[Normalized $\sigma_{12}$\label{subfig-5:test11}]{\ifrecompiletikz\tikzsetnextfilename{fig_06_c}\tikzexternalenable\input{images/fig_06_c}\tikzexternaldisable\else\includegraphics{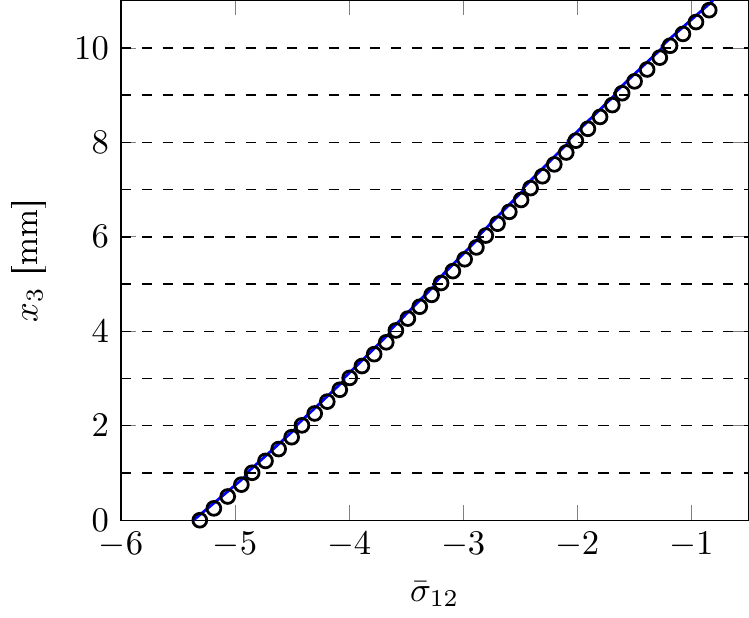}\fi}\hfill
	\subfigure[Normalized $\sigma_{33}$\label{subfig-6:test11}]{\ifrecompiletikz\tikzsetnextfilename{fig_06_f}\tikzexternalenable\input{images/fig_06_f}\tikzexternaldisable\else\includegraphics{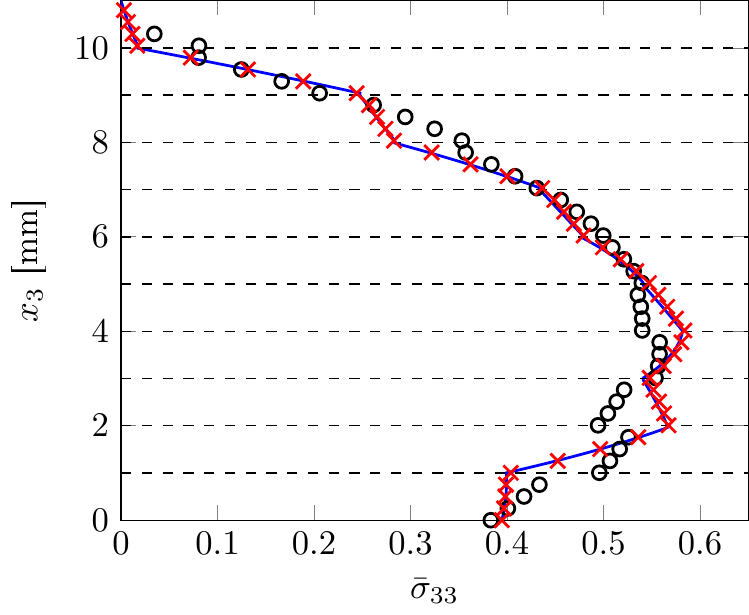}\fi}
	\caption{Through-the-thickness stress profiles evaluated at ($X_1=L/3$, $\theta=\bar{\theta}/3$) for IGA-Galerkin (degree of approximation $p=q=4$, $r=3$, and 22x22x4 control points). Case: Tube with mean radius-to-thickness ratio $S=20$, 11 layers, and $L=\bar{R}$ (\blueline~overkill IGA layerwise solution,~\blackcircles~single-element approach solution without post-processing, \redcrosses~post-processed solution).}
	\label{fig:LWI_IGAG_stress_pt1_1_l11_S20}
\end{figure}
\newpage
\subsubsection{Homogenized IGA-Collocation approach results}\label{subsec:homog_res}
\begin{figure}[!htbp]
	\centering
	\ifrecompiletikz\tikzsetnextfilename{fig_08}\tikzexternalenable\input{images/fig_08}\tikzexternaldisable\else\includegraphics{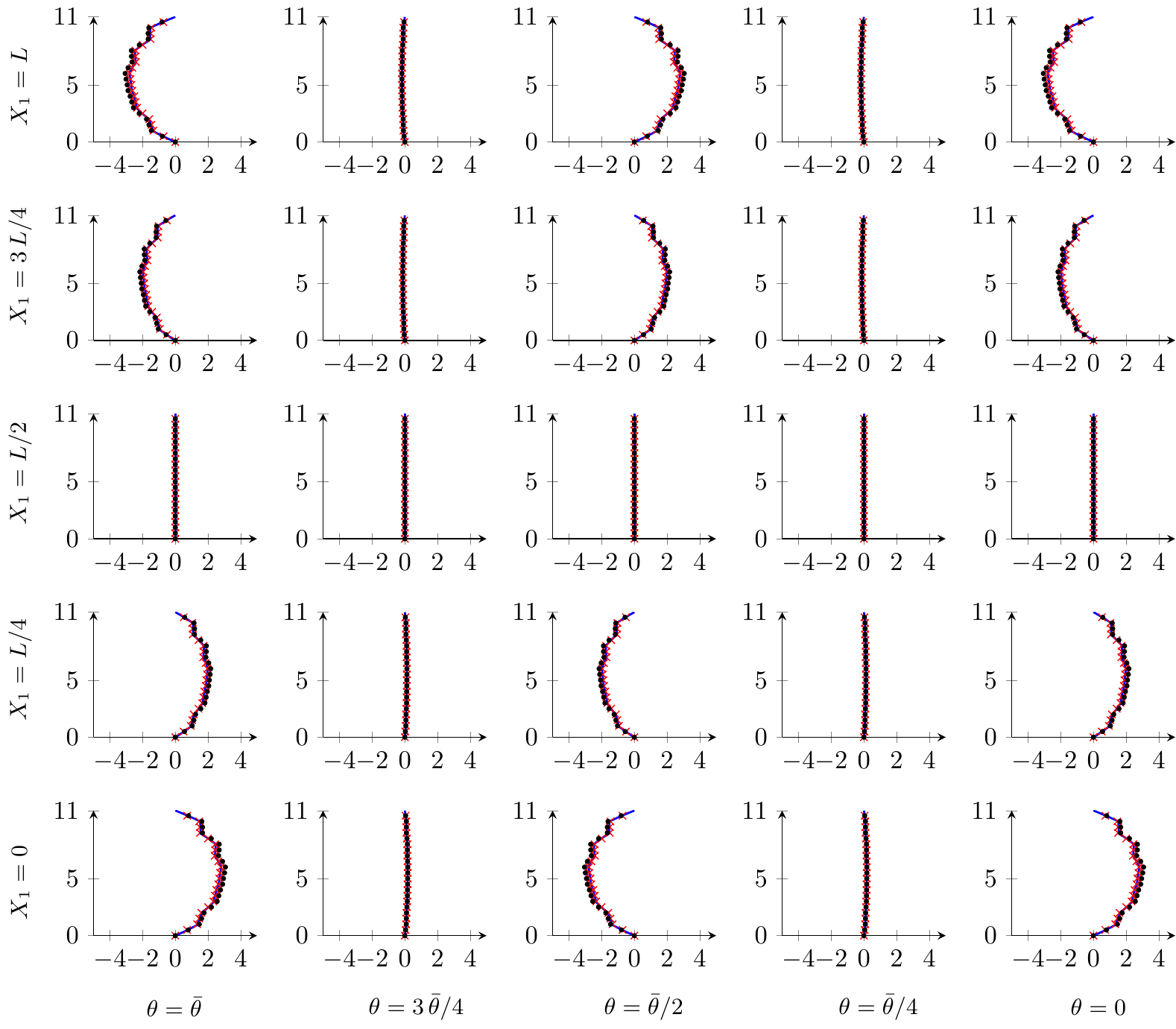}\fi
	\caption{Through-the-thickness $\bar{\sigma}_{13}$ profiles for several in plane sampling points ($S=\bar{R}/h=20$, $L=\bar{R}=220\,\text{mm}$, $h=11\,\text{mm}$), \redcrosses~IGA-Galerkin post-processed solution ($p=q=4$, $r=3$, and 22x22x4 control points),~\blackcirclesfull~IGA-Collocation  post-processed solution ($p=q=6$, $r=4$, and 22x22x5 control points),~\blueline~overkill IGA layerwise solution ($p=q=6$, $r=4$, and 36x36x55 control points).}
	\label{fig:LWI_IGAG_IGAC_samplingS13_l11_S20}
\end{figure}

\begin{figure}[!htbp]
	\centering
	\ifrecompiletikz\tikzsetnextfilename{fig_09}\tikzexternalenable\input{images/fig_09}\tikzexternaldisable\else\includegraphics{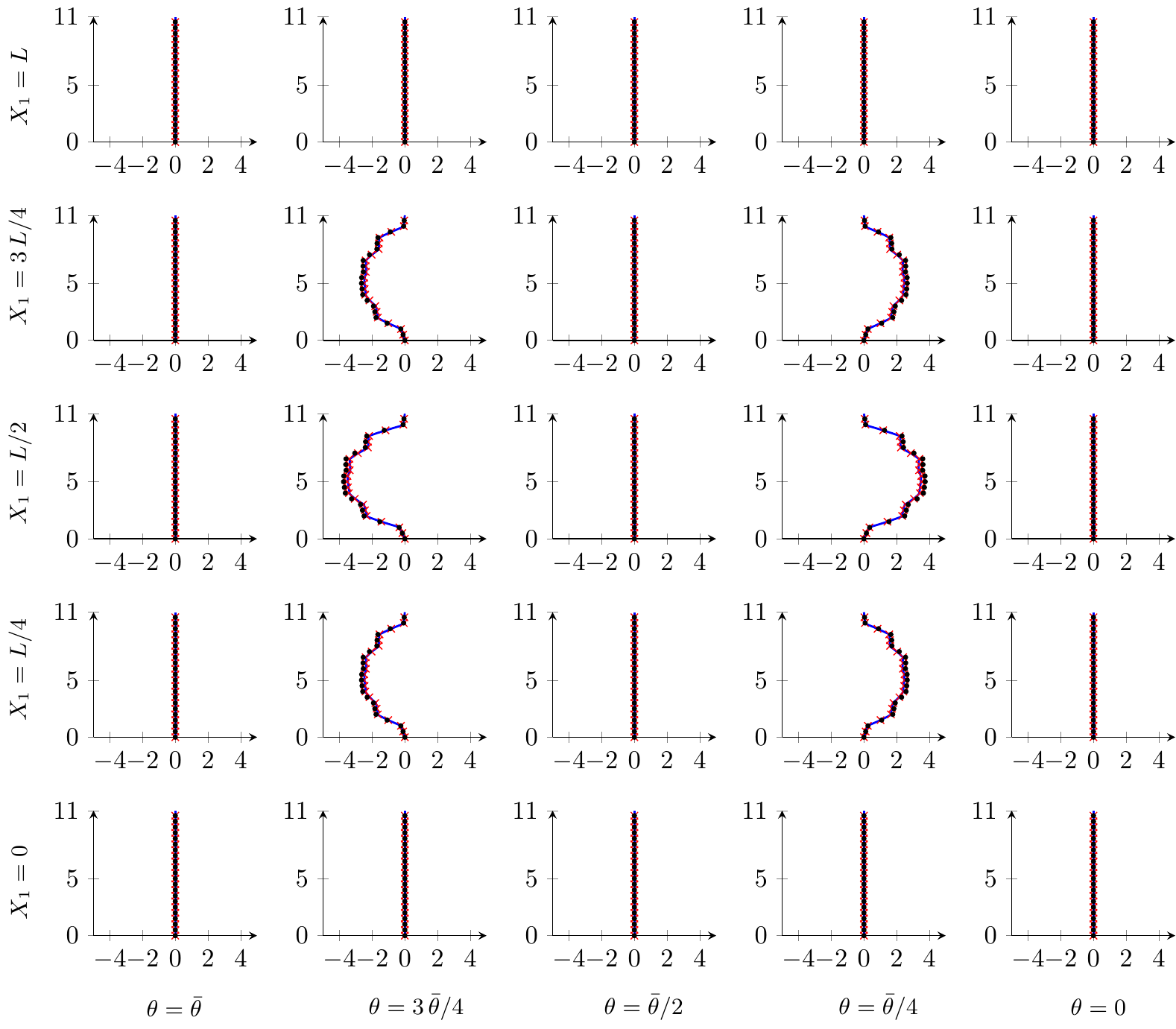}\fi
	\caption{Through-the-thickness $\bar{\sigma}_{23}$ profiles for several in plane sampling points ($S=\bar{R}/h=20$, $L=\bar{R}=220\,\text{mm}$, $h=11\,\text{mm}$), \redcrosses~IGA-Galerkin post-processed solution ($p=q=4$, $r=3$, and 22x22x4 control points),~\blackcirclesfull~IGA-Collocation  post-processed solution ($p=q=6$, $r=4$, and 22x22x5 control points),~\blueline~overkill IGA layerwise solution ($p=q=6$, $r=4$, and 36x36x55 control points).}
	\label{fig:LWI_IGAG_IGAC_samplingS23_l11_S20}
\end{figure}

\begin{figure}[!htbp]
	\centering
	\ifrecompiletikz\tikzsetnextfilename{fig_10}\tikzexternalenable\input{images/fig_10}\tikzexternaldisable\else\includegraphics{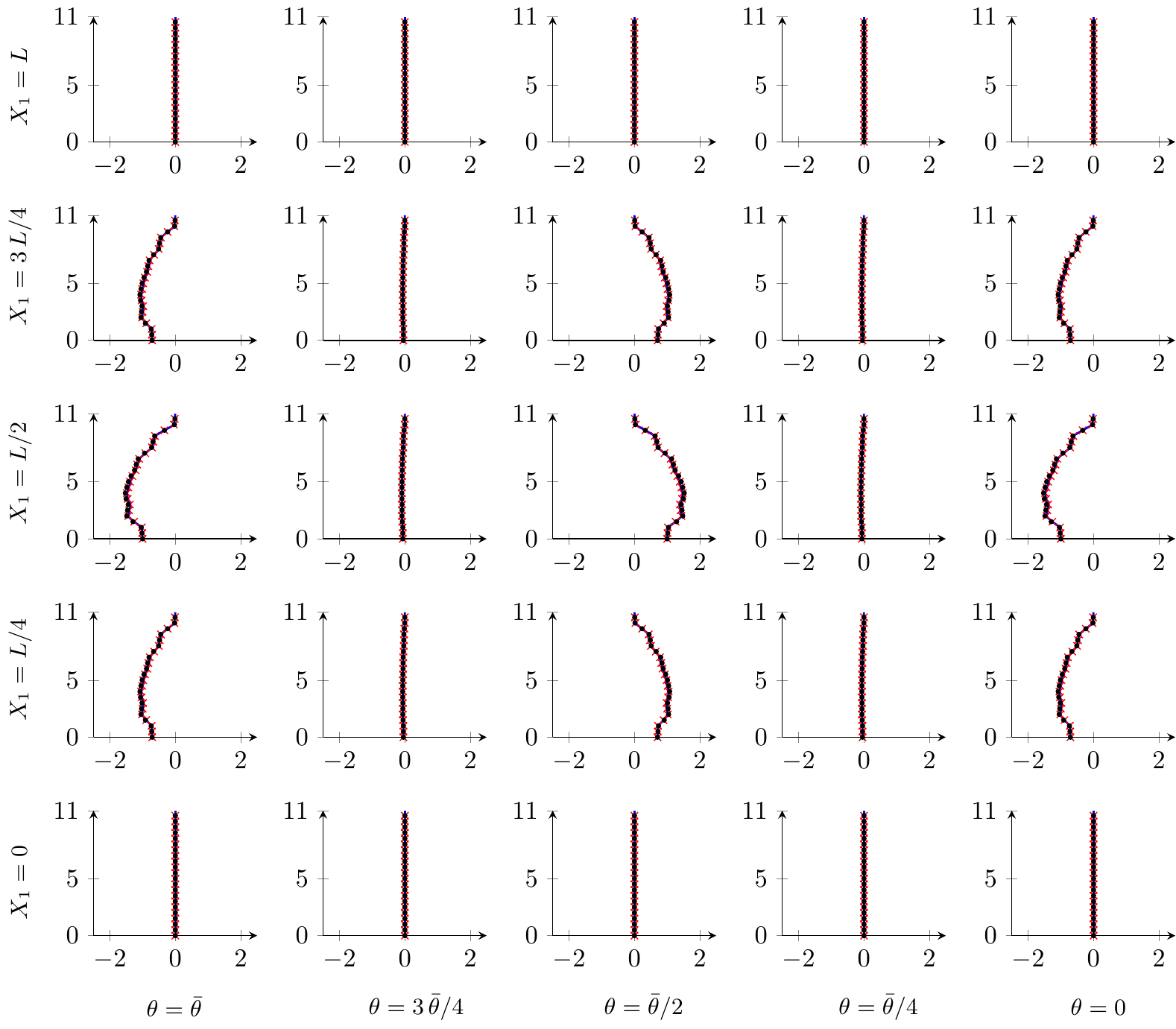}\fi
	\caption{Through-the-thickness $\bar{\sigma}_{33}$ profiles for several in plane sampling points ($S=\bar{R}/h=20$, $L=\bar{R}=220\,\text{mm}$, $h=11\,\text{mm}$), \redcrosses~IGA-Galerkin post-processed solution ($p=q=4$, $r=3$, and 22x22x4 control points),~\blackcirclesfull~IGA-Collocation  post-processed solution ($p=q=6$, $r=4$, and 22x22x5 control points),~\blueline~overkill IGA layerwise solution ($p=q=6$, $r=4$, and 36x36x55 control points).}
	\label{fig:LWI_IGAG_IGAC_samplingS33_l11_S20}
\end{figure}

In this section we report and discuss the results obtained, at the same sampling point ($X_1=L/3$, $\theta=\bar{\theta}/3$) as in Figure~\ref{fig:LWI_IGAG_stress_pt1_1_l11_S20}, from an IGA-Collocation homogenized displacement solution. Thus, in Figure~\ref{fig:LWI_IGAC_stress_pt1_1_l11_S20} we show in the left column the in-plane stress results which prove to be in good agreement with the reference IGA layerwise solution. In the right column, we depict out-of-plane stresses either obtained from the appropriate constitutive law or recovered from equilibrium equations. The overall solution, both in-plane and out-of-plane, is comparable to the reference one even though less precise than the one provided starting from IGA-Galerkin with the special integration rule. As further numerical tests exhibit the same kind of response provided by IGA-Collocation also for homogenized IGA-Galerkin, we believe that this behavior is most probably due to the material through-the-thickness homogenization strategy, which leads to a less accurate but still very good solution, especially for the $\sigma_{23}$ component which represent the preminent out-of-plane stress for this type of benchmark.  
\begin{figure}[!htbp]
	\centering
	\subfigure[Normalized $\sigma_{11}$\label{subfig-1:Colltest11S20}]{\ifrecompiletikz\tikzsetnextfilename{fig_07_a}\tikzexternalenable\input{images/fig_07_a}\tikzexternaldisable\else\includegraphics{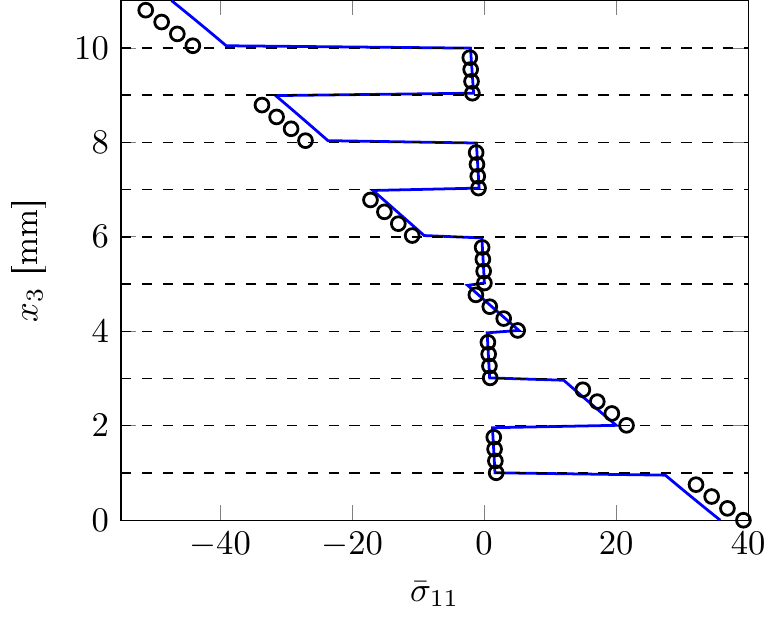}\fi}\hfill
	\subfigure[Normalized $\sigma_{13}$\label{subfig-2:Colltest11S20}]{\ifrecompiletikz\tikzsetnextfilename{fig_07_d}\tikzexternalenable\input{images/fig_07_d}\tikzexternaldisable\else\includegraphics{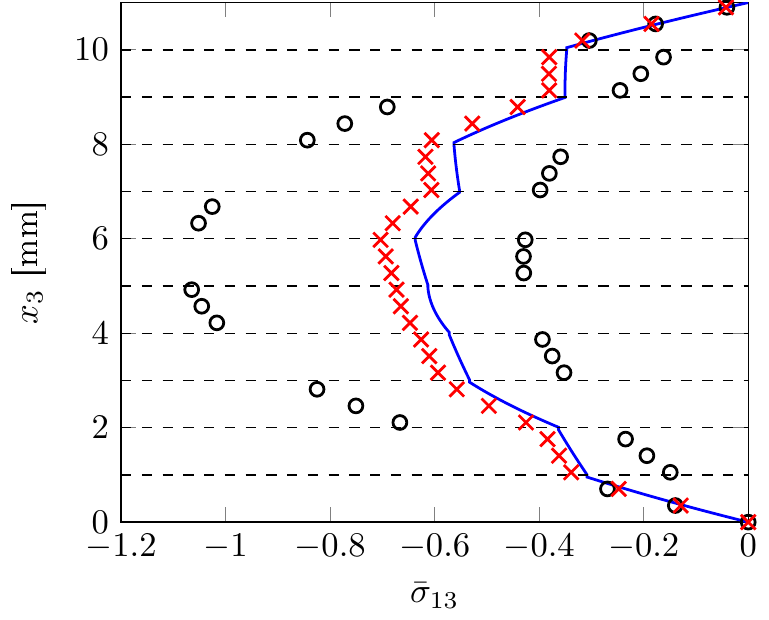}\fi}
	\caption{Through-the-thickness stress profiles evaluated at ($X_1=L/3$, $\theta=\bar{\theta}/3$) for IGA-Collocation (degree of approximation $p=q=6$, $r=4$, and 22x22x5 control points). Case: Tube with mean radius-to-thickness ratio $S=20$, 11 layers, and $L=\bar{R}$ (\blueline~overkill IGA layerwise solution,~\blackcircles~single-element approach solution without post-processing, \redcrosses~post-processed solution).}
	\end{figure}
	\begin{figure}[!htbp]
	\centering
	\ContinuedFloat
	\captionsetup{list=off,format=cont}
	\subfigure[Normalized $\sigma_{22}$\label{subfig-3:Colltest11S20}]{\ifrecompiletikz\tikzsetnextfilename{fig_07_b}\tikzexternalenable\input{images/fig_07_b}\tikzexternaldisable\else\includegraphics{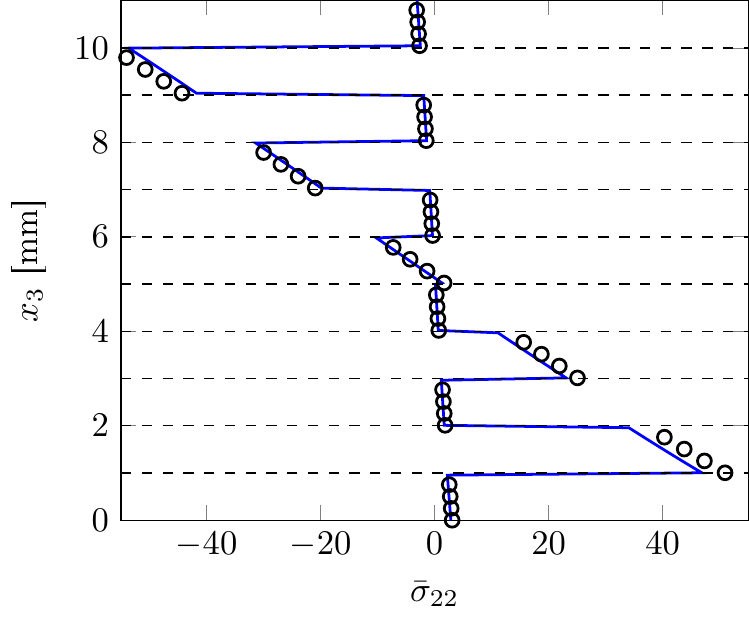}\fi}\hfill
	\subfigure[Normalized $\sigma_{23}$\label{subfig-4:Colltest11S20}]{\ifrecompiletikz\tikzsetnextfilename{fig_07_e}\tikzexternalenable\input{images/fig_07_e}\tikzexternaldisable\else\includegraphics{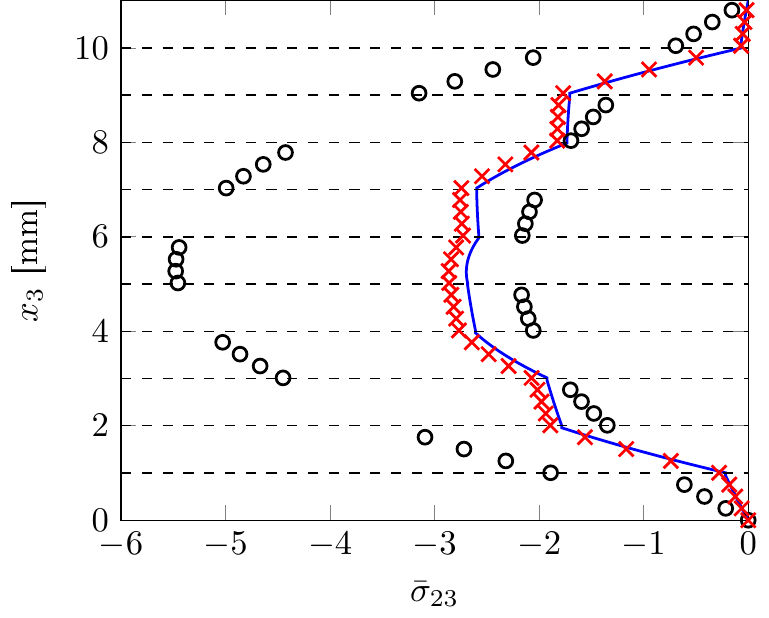}\fi}\\
	\subfigure[Normalized $\sigma_{12}$\label{subfig-5:Colltest11S20}]{\ifrecompiletikz\tikzsetnextfilename{fig_07_c}\tikzexternalenable\input{images/fig_07_c}\tikzexternaldisable\else\includegraphics{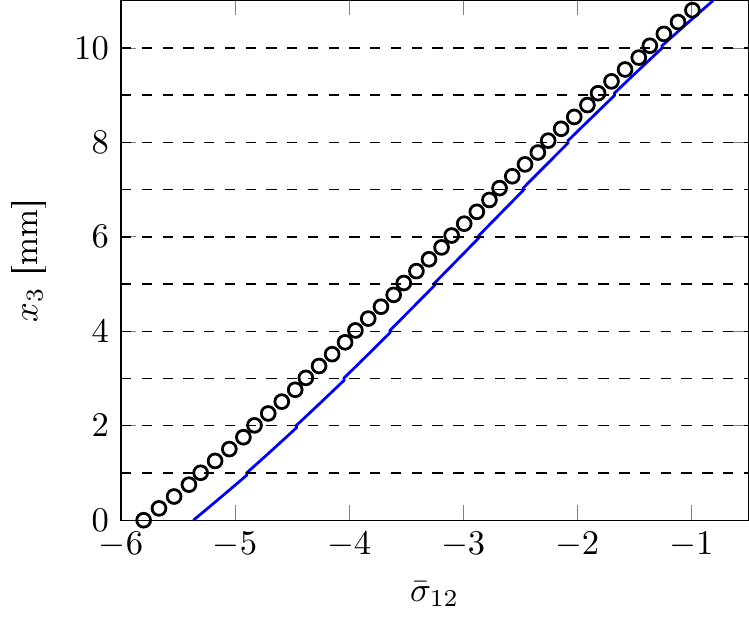}\fi}\hfill
	\subfigure[Normalized $\sigma_{33}$\label{subfig-6:Colltest11S20}]{\ifrecompiletikz\tikzsetnextfilename{fig_07_f}\tikzexternalenable\input{images/fig_07_f}\tikzexternaldisable\else\includegraphics{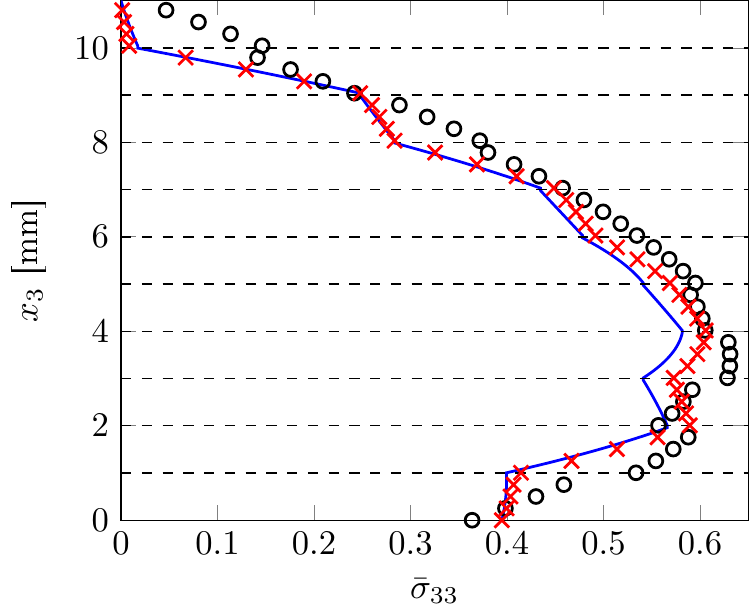}\fi}
	\caption{Through-the-thickness stress profiles evaluated at ($X_1=L/3$, $\theta=\bar{\theta}/3$) for IGA-Collocation (degree of approximation $p=q=6$, $r=4$, and 22x22x5 control points). Case: Tube with mean radius-to-thickness ratio $S=20$, 11 layers, and $L=\bar{R}$ (\blueline~overkill IGA layerwise solution,~\blackcircles~single-element approach solution without post-processing, \redcrosses~post-processed solution).}
	\label{fig:LWI_IGAC_stress_pt1_1_l11_S20}
\end{figure}
\newpage
In Figures~\ref{fig:LWI_IGAG_IGAC_samplingS13_l11_S20}-\ref{fig:LWI_IGAG_IGAC_samplingS33_l11_S20} we compare recovered out-of-plane stress profiles, either obtained via a displacement-based IGA-Galerkin or IGA-Collocation solution, against the considered reference layerwise results. To this end, the composite solid cylinder is sampled every quarter of length in both in-plane directions, showing the effectiveness of the presented post-processing strategy.
Also, this approach accurately captures the behavior of the cylinder at the boundaries, namely satisfies Neumann boundary conditions for transverse shear stresses at
the inner and outer surfaces of the laminate as in Equation~\eqref{eq:NBCs} and the symmetry conditions with respect to $\sigma_{13}$ and $\sigma_{23}$ as in Equation~\eqref{eq:sym_BCs}.
\newpage
\subsubsection{Parametric study on mean radius-to-thickness cylinder ratio}
The proposed approach is further tested considering different benchmarks with a varying slenderness parameter (\ie $S=20,\;30,\;40,\;\text{and}\;50$) for 11 and 33 layers and an increasing number of in-plane control points.

As a measure to assess the performance of the proposed stress recovery we introduce the following relative maximum error definition along the thickness
\begin{alignat}{2}\label{eq:error}
&\text{e}(\sigma_{i3})=\cfrac{\text{max}|\sigma_{i3}^\text{layerwise}(\bar{x}_1,\bar{x}_2,x_3)-\sigma_{i3}^\text{recovered}(\bar{x}_1,\bar{x}_2,x_3)|}{\text{max}|\sigma_{i3}^\text{layerwise}(\bar{x}_1,\bar{x}_2,x_3)|} \quad&&\quad i=1,2,3\,,
\end{alignat}
for which in the case of a zero stress profile we just compute the absolute maximum error to avoid 0/0 division. Also, in Equation~\eqref{eq:error} with the term $\sigma_{i3}^\text{layerwise}(\bar{x}_1,\bar{x}_2,x_3)$ ($i=1,2,3$) we mean that we take as a reference solution out-of-plane stress patterns obtained via an overkill layerwise approach which comprise again a degree of approximation $p=q=6$, $r=4$ and number of control points equal to 36x36x5x$n_l$.   

\begin{figure}[!htbp]
	\centering
	\hspace{-7pt}
	\subfigure[$\sigma_{13}$, 11 layers \label{subfig-1:Gal_errors13_l11_443_varS}]{\ifrecompiletikz\tikzsetnextfilename{fig_11_d}\tikzexternalenable\input{images/fig_11_d}\tikzexternaldisable\else\includegraphics{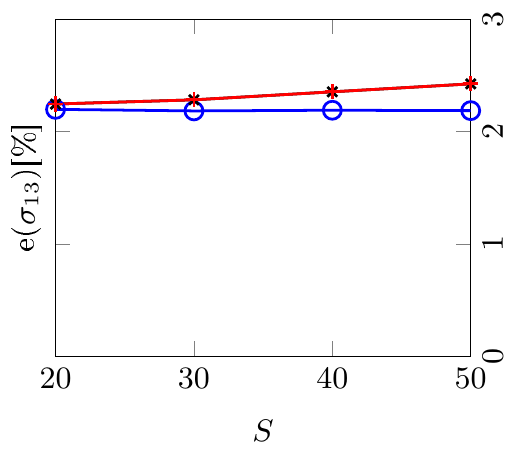}\fi}\hspace{-7pt}
	\subfigure[$\sigma_{23}$, 11 layers \label{subfig-2:Gal_errors23_l11_443_varS}]{\ifrecompiletikz\tikzsetnextfilename{fig_11_e}\tikzexternalenable\input{images/fig_11_e}\tikzexternaldisable\else\includegraphics{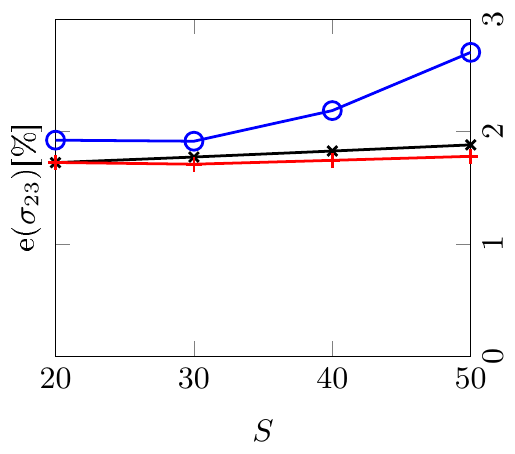}\fi}\hspace{-7pt}
	\subfigure[$\sigma_{33}$, 11 layers \label{subfig-3:Gal_errors33_l11_443_varS}]{\ifrecompiletikz\tikzsetnextfilename{fig_11_f}\tikzexternalenable\input{images/fig_11_f}\tikzexternaldisable\else\includegraphics{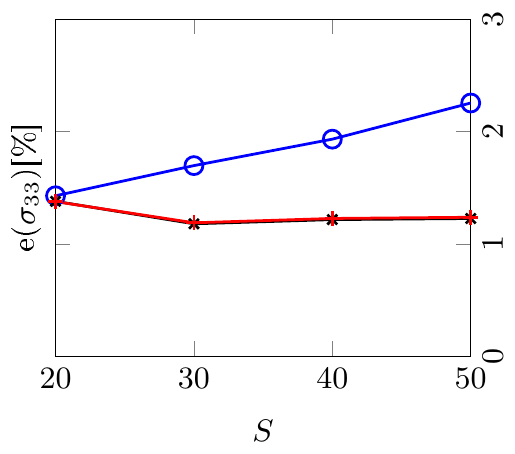}\fi}\\
	\hspace{-7pt}
	\subfigure[$\sigma_{13}$, 33 layers \label{subfig-4:Gal_errors13_l33_443_varS}]{\ifrecompiletikz\tikzsetnextfilename{fig_12_d}\tikzexternalenable\input{images/fig_12_d}\tikzexternaldisable\else\includegraphics{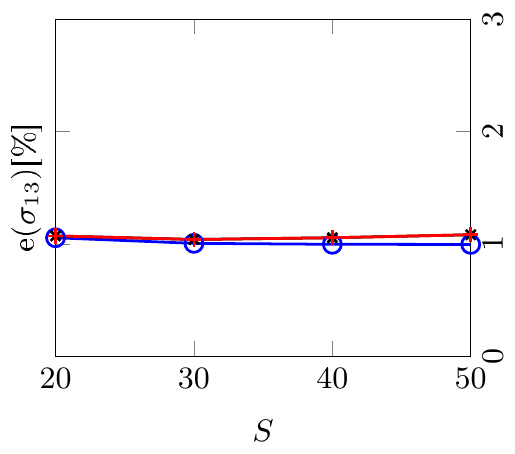}\fi}\hspace{-7pt}
	\subfigure[$\sigma_{23}$, 33 layers \label{subfig-5:Gal_errors23_l33_443_varS}]{\ifrecompiletikz\tikzsetnextfilename{fig_12_e}\tikzexternalenable\input{images/fig_12_e}\tikzexternaldisable\else\includegraphics{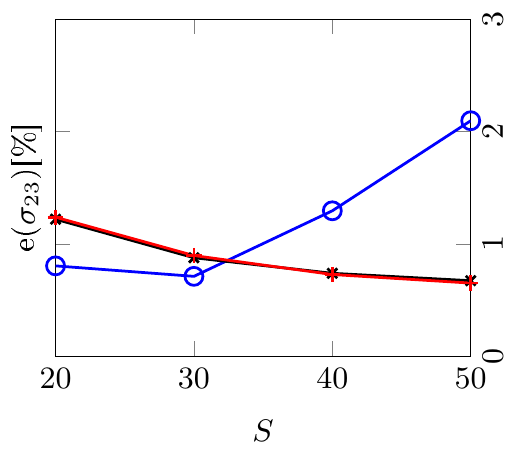}\fi}\hspace{-7pt}
	\subfigure[$\sigma_{33}$, 33 layers \label{subfig-6:Gal_errors33_l33_443_varS}]{\ifrecompiletikz\tikzsetnextfilename{fig_12_f}\tikzexternalenable\input{images/fig_12_f}\tikzexternaldisable\else\includegraphics{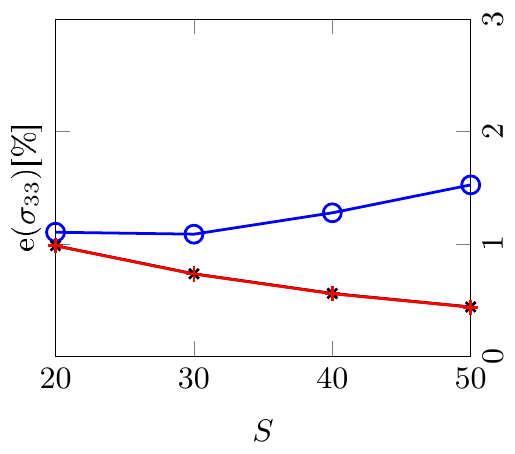}\fi}
    \caption{Maximum relative percentage error evaluation at ($X_1=L/3$, $\theta=\bar{\theta}/3$) of IGA-Galerkin post-processed single-element approach (degree of approximation $p=q=4$, $r=3$), with respect to an overkill IGA layerwise solution (degree of approximation $p=q=6$, $r=4$ and number of control points equal to 36x36x55). Different mean radius-to-thickness ratios $S$ are investigated for a number of layers equal to 11 and 33 (Number of in-plane control points per parametric direction:~\bluesolidcircle~11,~\blacksolidx~22,~\redsolidcross~44). }
	\label{fig:Gal_conv443_varS}
\end{figure}

In Figure~\ref{fig:Gal_conv443_varS} we assess the performance of the isogeometric Galerkin approach coupled with the presented post-processing technique at $(X_1=L/3, \theta=\bar{\theta}/3)$, while in Figure~\ref{fig:Col_conv666_varS} we test the stress recovery technique applied to the described homogenized IGA collocation method for the same in-plane sampling point.

The post-processing approach seems to be particularly
suitable to tackle the behavior of slender laminates characterized by a significant number of layers, leading to errors that are typically in the 10\% range or lower at convergence despite the chosen displacement-based approach. More specifically IGA-Galekin provides errors in the order of 2\% and 1\% in average for the analyzed 11 and 33 layer cases, respectively, while IGA-collocation allows to obtain errors in the range of 5\% or lower, whenever $S\ge20$.
Furthermore, collocation proves to be more sensitive to variations of the slenderness parameter $S$ with respect to the proposed Galerkin approach. This is due to the fact that a slender laminate made of a significant number of layers is closer to a structure with average material properties.
We would like to remark that the degrees of approximation adopted for the considered example are taken accordingly to \cite{Patton2019,Dufour2018}, which proved to be suitable choices to model less complex geometries.

Moreover, we observe that considering a mesh comprising 22x22 in-plane control points allows to correctly describe the considered cases and further in-plane refinement operations do not seem to significantly improve the quality of the overall solution. This is due to the fact that the modeling error, given by the equilibrium-based direct integration, dominates over the approximation one.
\begin{figure}[!htbp]
	\centering
	\hspace{-7pt}
	\subfigure[$\sigma_{13}$, 11 layers \label{subfig-1:Col_errors13_l11_666_varS}]{\ifrecompiletikz\tikzsetnextfilename{fig_13_d}\tikzexternalenable\input{images/fig_13_d}\tikzexternaldisable\else\includegraphics{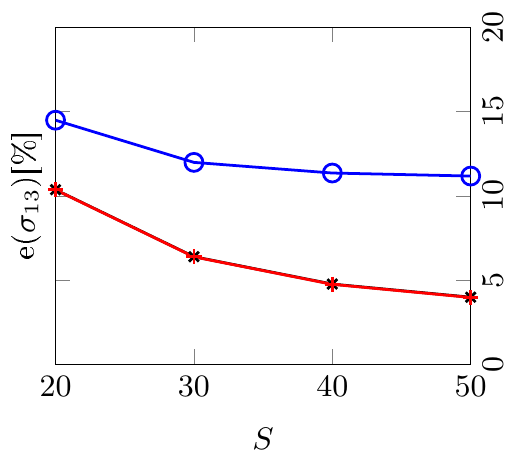}\fi}\hspace{-7pt}
	\subfigure[$\sigma_{23}$, 11 layers \label{subfig-2:Col_errors23_l11_666_varS}]{\ifrecompiletikz\tikzsetnextfilename{fig_13_e}\tikzexternalenable\input{images/fig_13_e}\tikzexternaldisable\else\includegraphics{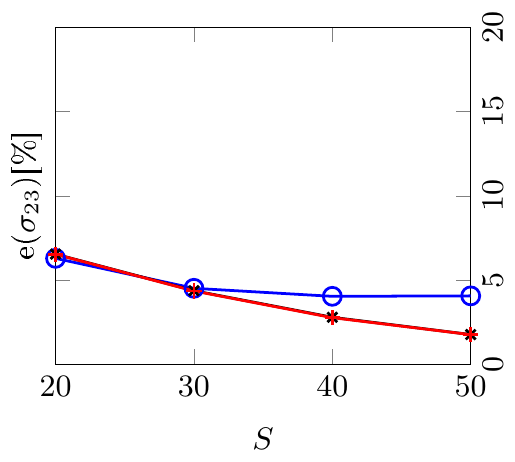}\fi}\hspace{-7pt}
	\subfigure[$\sigma_{33}$, 11 layers \label{subfig-3:Col_errors33_l11_666_varS}]{\ifrecompiletikz\tikzsetnextfilename{fig_13_f}\tikzexternalenable\input{images/fig_13_f}\tikzexternaldisable\else\includegraphics{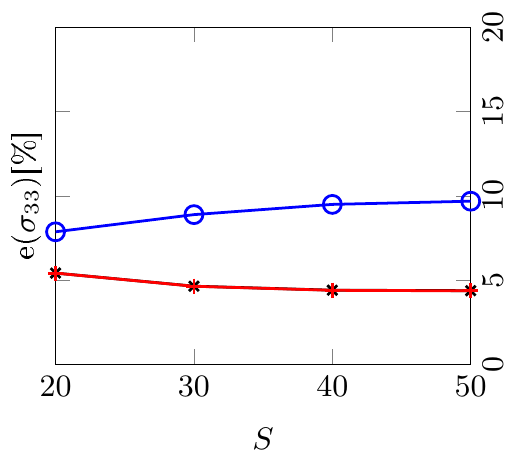}\fi}
	\caption{Maximum relative percentage error evaluation at ($X_1=L/3$, $\theta=\bar{\theta}/3$) of IGA-Collocation post-processed single-element approach (degree of approximation $p=q=4$, $r=3$), with respect to an overkill IGA layerwise solution (degree of approximation $p=q=6$, $r=4$ and number of control points equal to 36x36x55). Different mean radius-to-thickness ratios $S$ are investigated for a number of layers equal to 11 and 33 (Number of in-plane control points per parametric direction:~\bluesolidcircle~11,~\blacksolidx~22,~\redsolidcross~44).}
\end{figure}
\begin{figure}[!htbp]
	\centering
	\ContinuedFloat
	\captionsetup{list=off,format=cont}
	\hspace{-7pt}
	\subfigure[$\sigma_{13}$, 33 layers \label{subfig-4:Col_errors13_l33_666_varS}]{\ifrecompiletikz\tikzsetnextfilename{fig_14_d}\tikzexternalenable\input{images/fig_14_d}\tikzexternaldisable\else\includegraphics{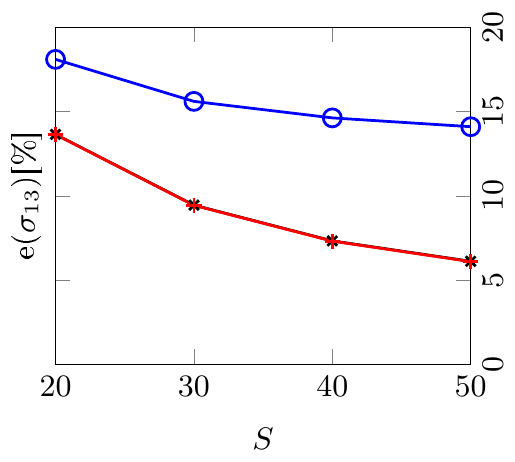}\fi}\hspace{-7pt}
	\subfigure[$\sigma_{23}$, 33 layers \label{subfig-5:Col_errors23_l33_666_varS}]{\ifrecompiletikz\tikzsetnextfilename{fig_14_e}\tikzexternalenable\input{images/fig_14_e}\tikzexternaldisable\else\includegraphics{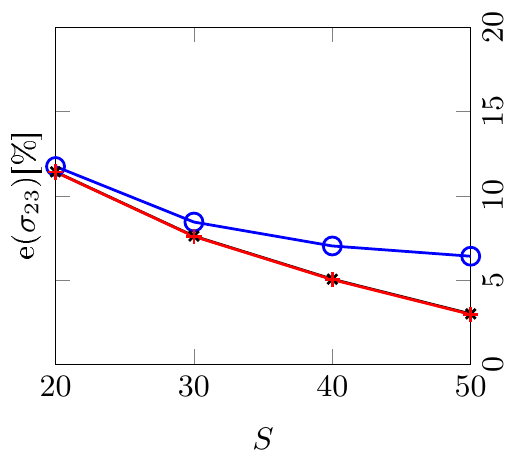}\fi}\hspace{-7pt}
	\subfigure[$\sigma_{33}$, 33 layers \label{subfig-6:Col_errors33_l33_666_varS}]{\ifrecompiletikz\tikzsetnextfilename{fig_14_f}\tikzexternalenable\input{images/fig_14_f}\tikzexternaldisable\else\includegraphics{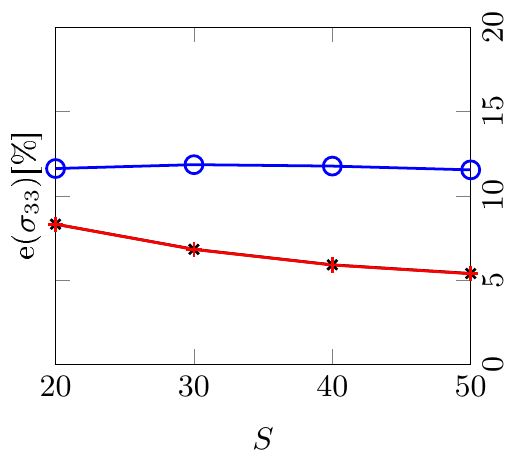}\fi}
	\caption{Maximum relative percentage error evaluation at ($X_1=L/3$, $\theta=\bar{\theta}/3$) of IGA-Collocation post-processed single-element approach (degree of approximation $p=q=4$, $r=3$), with respect to an overkill IGA layerwise solution (degree of approximation $p=q=6$, $r=4$ and number of control points equal to 36x36x55). Different mean radius-to-thickness ratios $S$ are investigated for a number of layers equal to 11 and 33 (Number of in-plane control points per parametric direction:~\bluesolidcircle~11,~\blacksolidx~22,~\redsolidcross~44).}
	\label{fig:Col_conv666_varS}
\end{figure}

Finally, to highlight the quantitative gain given by the proposed stress-recovery application, we report in Table~\ref{tab:errorl11} the relative maximum error at ($L/3$, $\bar{\theta}/3$) obtained considering out-of-plane stress components either directly computed via the appropriate constitutive law or reconstructed.
\begin{table}[!htbp]
	\caption{Simply supported solid composite cylinder under sinusoidal load with a number of layers equal to 11. Out-of-plane stress maximum relative error along the thickness with respect to an overkill IGA layerwise solution (degree of approximation $p=q=6$, $r=4$ and number of control points equal to 36x36x55) at ($L/3$, $\bar{\theta}/3$). Assessment of the proposed isogeometric collocation-based homogenized single element approach (IGA-Collocation) and the isogeometric Galerkin single element method (IGA-Galerkin) before and after the proposed post-processing technique application for different approximation degrees and a fixed mesh comprising 22x22 in-plane control points and a single element through the thickness.}
	\vspace{0.5cm}
	\centering	
	\begin{tabular}{? l | c ? c | c | c ? c | c | c ?}
		\thickhline
		\multicolumn{2}{? c ?}{Method}
		& \multicolumn{3}{c ?}{IGA-Galerkin}&\multicolumn{3}{c ?}{IGA-Collocation}\Tstrut\Bstrut\\
		\multicolumn{2}{? c ?}{Degree of approximation}&\multicolumn{3}{c ?}{ $p=q=4$, $r=3$}&\multicolumn{3}{c ?}{ $p=q=6$, $r=4$}\Tstrut\Bstrut\\\thickhline
		\multirow{2}{*}{$\boldsymbol{S}$} &\textbf{Out-of-plane maximum relative error}
		&e($\sigma_{13}$)
		&e($\sigma_{23}$) 
		&e($\sigma_{33}$) 
		&e($\sigma_{13}$) 
		&e($\sigma_{23}$) 
		&e($\sigma_{33}$)  
		\Tstrut\Bstrut\\\cline{3-8}
		& \textbf{at ($\boldsymbol{L/3}$, $\boldsymbol{\bar{\theta}/3}$, $\boldsymbol{x_3}$)} &[\%]&[\%]&[\%]&[\%]&[\%]&[\%]
		\Tstrut\Bstrut\\\thickhline
		\multirow{2}{*}{20} & before post-processing & 46.6 & 47.8 & 18.7 & 71.2 & 105 & 23.1\Tstrut\Bstrut\\
		& after post-processing & 2.25 & 1.72 & 1.38 & 10.4 & 6.54 & 5.41\Tstrut\Bstrut\\\hline
		\multirow{2}{*}{30} & before post-processing & 46.8 & 47.6 & 19.8 & 67.0 & 103 & 27.4\Tstrut\Bstrut\\
		& after post-processing & 2.28 & 1.77 & 1.18 & 6.39 & 4.37 & 4.63\Tstrut\Bstrut\\\hline
		\multirow{2}{*}{40} & before post-processing & 47.0 & 47.7 & 22.6 & 65.4 & 100 & 30.5\Tstrut\Bstrut\\
		& after post-processing & 2.36 & 1.83 & 1.22 & 4.76 & 2.79 & 4.40\Tstrut\Bstrut\\\hline
		\multirow{2}{*}{50} & before post-processing& 47.3 & 47.8 & 25.2 & 64.9 & 97.9 & 33.1\Tstrut\Bstrut\\
		& after post-processing & 2.43 & 1.88 & 1.23 & 3.98 & 1.76 & 4.37\Tstrut\Bstrut\\\thickhline
	\end{tabular}\label{tab:errorl11}
\end{table}

	\newpage
	\section{Conclusions}\label{sec:conclusions}
    In this work we propose an efficient two-step single-element displacement-based approach to model the stress field in laminated curved structures. This technique combines coarse 3D isogeometric computations, either using a calibrated layerwise integration rule or a homogenized approach, with an a-posteriori step where equilibrium is directly imposed in strong form. Both IGA Galerkin and collocation methods are successfully investigated. Homogenization is immediately effective only for symmetric ply distributions, which cover most common cases in practice. While this approach represents a natural choice for collocation, it can be regarded as a less accurate but cost-effective alternative for the proposed IGA Galerkin method, requiring only $r+1$ quadrature points through the thickness. In any case, despite being more computationally demanding, the proposed layerwise integration rule (i.e., $r+1$ Gauss points per ply) for the IGA Galerkin approach allows to correctly capture the behavior of composites for any stacking sequence and investigate more complex constitutive models.
    	
    Using only one element through the thickness with highly continuous shape functions provides a good in-plane solution, while the out-of-plane stress profiles violate the continuity requirements prescribed by equilibrium. Adopting a local system, which grants that no additional terms appear in the balance of linear momentum equations, an accurate solution is attainable also in terms of out-of-plane components after the post-processing step.  
    
    The proposed stress recovery results in a straightforward approach in terms of implementation since it is only based on the numerical integration through the thickness of equilibrium equations and all the required components can be easily computed differentiating the obtained displacement solution.
    The a-posteriori step requires the shape functions to be highly continuous (as is fully granted by IGA), is direct, and can be easily performed at locations of interest, resulting in a very interesting alternative to layerwise techniques especially in the case of slender composites made of a significant number of layers. Several benchmarks assess the behavior of the proposed approach also on the laminate boundary and its sensitivity to different parameters such as the mean radius-to-thickness ratio and number of layers.
     
    Further research topics currently under investigation consist of the extension of this approach to bivariate shells and the inclusion of large deformations. Also, more efficient through-the-thickness integration strategies will be considered in future studies.

	\section*{Acknowledgments}
	This work was partially supported by Ministero dell'Istruzione, dell'Universit\`a e della Ricerca through the project ``XFAST-SIMS: Extra fast and accurate simulation of complex structural systems'', within the program Progetti di ricerca di Rilevante Interesse Nazionale (PRIN). P. Antol\'in was also partially supported by the European Research Council through the H2020 ERC Advanced Grant 2015 n. 694515 CHANGE, and by the Swiss National Science Foundation through the project ``Design-through-Analysis (of PDEs): the litmus test'' n. 40B2-0 187094 (BRIDGE Discovery 2019). J. Kiendl was partially supported by the European Research Council through the H2020 ERC Consolidator Grant 2019 n. 864482 FDM$^2$. Finally, the authors would like to thank Dr. John-Eric Dufour (University of Texas at Arlington) for his support and comments on this paper.

	\appendix
	\section{Components of stress derivatives}\label{sec:stress_comp}
\subsection{Stress derivatives respect to the local reference system}
In Section \ref{sec:recovery}, the stress recovery procedure proposed in \cite{Dufour2018} is extended to the case of curved laminated structures.
In this regard, the stress tensor first and second derivatives in the chosen local reference system $\{x_\alpha\}$ need to be computed.
Therefore, starting from the obtained displacement field in the global cartesian system $\{X_i\}$ as described in Section~\ref{sec:numerical_strategies}, $\diff\stress/\diff x_\mu$, and $\diff^2\stress/\diff x_\mu\diff x_\nu$ have to be assessed.

Thus, from Equation~\eqref{eq:stress_strain}, the stress derivatives can be computed as
\begin{align}
\frac{\diff \stress}{\diff x_\mu} &= \frac{\diff \CC}{\diff x_\mu}:\strain +
\CC:\frac{\diff \strain}{\diff x_\mu}\,,\label{eq:derstress_local_expr}\\
\frac{\diff^2 \stress}{\diff x_\mu\diff x_\nu} &= \frac{\diff^2 \CC}{\diff x_\mu\diff x_\nu}:\strain
+ \frac{\diff \CC}{\diff x_\mu}:\frac{\diff \strain}{\diff x_\nu}
+ \frac{\diff \CC}{\diff x_\nu}:\frac{\diff \strain}{\diff x_\mu}
+ \CC:\frac{\diff^2 \strain}{\diff x_\mu\diff x_\nu}\,,\label{eq:der2stress_local_expr}
\end{align}
which need to be further detailed in terms of explicit expressions for the first and second derivatives of $\CC$ and $\strain$.

The strain and its derivatives can be easily obtained referred to the local system $\{x_\alpha\}$ as:
\begin{subequations}\label{eq:strain_derstrain_local}
\begin{align}
\strain &= \varepsilon_{\alpha\beta}\e_\alpha\otimes\e_\beta\,,\label{eq:strain_local}\\
\frac{\diff\strain}{\diff x_\mu} &= \varepsilon_{\alpha\beta,\mu}\e_\alpha\otimes\e_\beta\,,\label{eq:derstrain_local}\\
\frac{\diff^2\strain}{\diff x_\mu\diff x_\nu} &= \varepsilon_{\alpha\beta,\mu\nu}\e_\alpha\otimes\e_\beta\,.\label{eq:der2strain_local}
\end{align}
\end{subequations}
From Eq.~\eqref{eq:strain_components}, and applying the chain rule, the components $\varepsilon_{\alpha\beta}$, $\varepsilon_{\alpha\beta,\mu}$, and $\varepsilon_{\alpha\beta,\mu\nu}$ are computed as:
\begin{subequations}\label{eq:strain_derstrain_local_comp}
\begin{align}
\varepsilon_{\alpha\beta} &= \tilde{\varepsilon}_{ij}C_{i\alpha} C_{j\beta} = \frac{1}{2}\left(\tilde u_{i,j} + \tilde u_{j,i}\right) C_{i\alpha} C_{j\beta}\,,\label{eq:strain_local_comp}\\
\varepsilon_{\alpha\beta,\mu} &= \frac{\diff\varepsilon_{\alpha\beta}}{\diff x_\mu} =
\frac{\diff\varepsilon_{\alpha\beta}}{\diff X_k}\frac{\diff X_k}{\diff x_\mu} =
\frac{1}{2}\left(\tilde u_{i,jk} + \tilde u_{j,ik}\right) C_{i\alpha} C_{j\beta} C_{k\mu}\,,\label{eq:derstrain_local_comp}\\
\varepsilon_{\alpha\beta,\mu\nu} &= \frac{\diff^2\varepsilon_{\alpha\beta}}{\diff x_\mu\diff x_\nu} =
\frac{\varepsilon_{\alpha\beta,\mu}}{\diff X_l}\frac{\diff X_l}{\diff x_\nu} =
\frac{1}{2}\left(\tilde u_{i,jkl} + \tilde u_{j,ikl}\right) C_{i\alpha} C_{j\beta} C_{k\mu} C_{l\nu}\,,\label{eq:der2strain_local_comp}
\end{align}
\end{subequations}
where we applied the definition \eqref{eq:metric} and the property \eqref{eq:metric_ders}.
The displacement derivatives $\tilde u_{i,j}$, $\tilde u_{i,jk}$, and $\tilde u_{i,jkl}$ are simply computed as:
\begin{subequations}
\begin{align}
\tilde u_{i,j} &= \frac{\diff \tilde u_{i}}{\diff X_j}\,,\label{eq:der1displ_local_comp}\\
\tilde u_{i,jk} &= \frac{\diff^2 \tilde u_{i}}{\diff X_j\diff X_k}\,,\label{eq:der2displ_local_comp}\\
\tilde u_{i,jkl} &= \frac{\diff^3 \tilde u_{i}}{\diff X_j\diff X_k\diff X_l}\,.\label{eq:der3displ_local_comp}
\end{align}
\end{subequations}
Furthermore, we would like to remark that the derivatives $\tilde{u}_{i,j}$, $\tilde{u}_{i,jk}$, and $\tilde{u}_{i,jkl}$ are evaluated starting from the global components of the solution $\tilde{u}_i$ in a straightforward manner. In fact, relying on high order continuity properties, IGA allows to easily perform the necessary computations.

In the same way, the derivatives of $\CC$, defined starting from Equation~\eqref{eq:CC_local}, are computed as
\begin{subequations}\label{eq:CC_der_local}
\begin{align}
\frac{\diff \CC}{\diff x_\mu} &=
  \CC_{\alpha\beta\gamma\delta,\mu} \a_\alpha\otimes\a_\beta\otimes\a_\gamma\otimes\a_\delta +
  \CC_{\alpha\beta\gamma\delta}\frac{\diff}{\diff x_\mu}\left(\a_\alpha\otimes\a_\beta\otimes\a_\gamma\otimes\a_\delta\right),\label{eq:derCC_local}\\
\begin{split}
\frac{\diff^2 \CC}{\diff x_\mu\diff x_\nu} &=
  \CC_{\alpha\beta\gamma\delta,\mu\nu} \a_\alpha\otimes\a_\beta\otimes\a_\gamma\otimes\a_\delta
  + \CC_{\alpha\beta\gamma\delta,\mu} \frac{\diff}{\diff x_\nu}\left(\a_\alpha\otimes\a_\beta\otimes\a_\gamma\otimes\a_\delta\right)\\
  + &\CC_{\alpha\beta\gamma\delta,\nu} \frac{\diff}{\diff x_\mu}\left(\a_\alpha\otimes\a_\beta\otimes\a_\gamma\otimes\a_\delta\right)
  + \CC_{\alpha\beta\gamma\delta}\frac{\diff^2}{\diff x_\mu\diff x_\nu}\left(\a_\alpha\otimes\a_\beta\otimes\a_\gamma\otimes\a_\delta\right)\,,
\end{split}\label{eq:der2CC_local}
\end{align}
\end{subequations}
where
\begin{subequations}
\begin{align}
\CC_{\alpha\beta\gamma\delta,\mu} &= \frac{\diff\CC_{\alpha\beta\gamma\delta}}{\diff x_\mu}\,,\\
\CC_{\alpha\beta\gamma\delta,\mu\nu} &= \frac{\diff^2\CC_{\alpha\beta\gamma\delta}}{\diff x_\mu\diff x_\nu}\,.
\end{align}
\end{subequations}
In Eq.~\eqref{eq:CC_der_local}, the derivatives $\diff\a_\alpha/\diff x_\mu$ and $\diff^2\a_\alpha/\diff x_\mu\diff x_\nu$ do not vanish in general, as already discussed in Section \ref{subsec:kin}.
And for heterogeneous materials, in which the material coefficients may change from point to point, the terms $\CC_{\alpha\beta\gamma\delta,\mu}$ and $\CC_{\alpha\beta\gamma\delta,\mu\nu}$ may be also different from zero.
\begin{remark}\label{remark:5}
For the case of homogeneous anisotropic materials (including orthotropic ones)
in  which the fibers have a constant orientation with 
respect to the local basis $\{\a_\alpha\}$,
the terms $\CC_{\alpha\beta\gamma\delta,\mu}$ and $\CC_{\alpha\beta\gamma\delta,\mu\nu}$ vanish.
This is the case, e.g., of multi-layered structures in which stacks of materials with different orientations are used: Within each layer $(k)$ the material coefficients $\CC^{(k)}_{\alpha\beta\gamma\delta}$ are constant for the chosen basis $\{\a_\alpha\}$.
\end{remark}
Finally, the stress derivative terms $\sigma_{\alpha\beta,\mu}$ and $\sigma_{\alpha\beta,\mu\nu}$ required by the recovery integrals \eqref{eq:recoveryShear} and \eqref{eq:s33,3} are computed as:
\begin{subequations}\label{eq:stress_der_comp}
\begin{align}
\sigma_{\alpha\beta,\mu} &=\left(\e_\alpha\otimes\e_\beta\right):\frac{\diff \CC}{\diff x_\mu}:\strain +
\left(\e_\alpha\otimes\e_\beta\right) : \CC:\frac{\diff \strain}{\diff x_\mu}\,,\\
\begin{split}
\sigma_{\alpha\beta,\mu\nu} &=  \left(\e_\alpha\otimes\e_\beta\right) :\frac{\diff^2 \CC}{\diff x_\mu\diff x_\nu}:\strain
+ \left(\e_\alpha\otimes\e_\beta\right) :\frac{\diff \CC}{\diff x_\mu}:\frac{\diff \strain}{\diff x_\nu}\\
&+ \left(\e_\alpha\otimes\e_\beta\right) :\frac{\diff \CC}{\diff x_\nu}:\frac{\diff \strain}{\diff x_\mu}
+ \left(\e_\alpha\otimes\e_\beta\right) : \CC:\frac{\diff^2 \strain}{\diff x_\mu\diff x_\nu}\,.
\end{split}
\end{align}
\end{subequations}
Thus, considering the Eqs.~\eqref{eq:CC_der_local}, the expression terms $\sigma_{\alpha\beta,\mu}$ and $\sigma_{\alpha\beta,\mu\nu}$ can be written as:
\begin{subequations}\label{eq:stress_der_local_comp}
\begin{align}
\begin{split}
\sigma_{\alpha\beta,\mu} &= \Big(
\CC_{\alpha\beta\gamma\delta,\mu}
+ \CC_{\psi\beta\gamma\delta}A_{\psi \alpha\mu}
+ \CC_{\alpha\psi \gamma\delta}A_{\psi \beta\mu}
+ \CC_{\alpha\beta\psi\delta}A_{\psi \gamma\mu}\\
&+ \CC_{\alpha\beta\gamma\psi}A_{\psi \delta\mu}
\Big)\varepsilon_{\gamma\delta}
+\CC_{\alpha\beta\gamma\delta} \varepsilon_{\gamma\delta,\mu}\,,
\end{split}\label{eq:derstress_local_comp}\\
\begin{split}
\sigma_{\alpha\beta,\mu\nu} &= \Big(
\CC_{\alpha\beta\gamma\delta,\mu\nu}\\
&+ \CC_{\psi \beta\gamma\delta,\mu}A_{\psi \alpha\nu}
+ \CC_{\alpha\psi \gamma\delta,\mu}A_{\psi \beta\nu}
+ \CC_{\alpha\beta\psi\delta,\mu}A_{\psi \gamma\nu}
+ \CC_{\alpha\beta\gamma\psi,\mu}A_{\psi \delta\nu}\\
&+ \CC_{\psi\beta\gamma\delta,\nu}A_{\psi \alpha\mu}
+ \CC_{\alpha\psi\gamma\delta,\nu}A_{\psi \beta\mu}
+ \CC_{\alpha\beta\psi\delta,\nu}A_{\psi \gamma\mu}
+ \CC_{\alpha\beta\gamma\psi,\nu}A_{\psi \delta\mu}\\
&+ \CC_{\psi\beta\gamma\delta}B_{\psi \alpha\mu\nu}
+ \CC_{\psi\omega\gamma\delta}A_{\psi \alpha\mu}A_{\omega \beta\nu}
+ \CC_{\psi\beta\omega\delta}A_{\psi \alpha\mu}A_{\omega \gamma\nu}
+ \CC_{\psi\beta\gamma\omega}A_{\psi \alpha\mu}A_{\omega \delta\nu}\\
&+ \CC_{\omega\psi\gamma\delta}A_{\psi \beta\mu}A_{\omega \alpha\nu}
+ \CC_{\alpha\psi\gamma\delta}B_{\psi \beta\mu\nu}
+ \CC_{\delta\mu\psi\omega\delta}A_{\psi \beta\mu}A_{\omega \gamma\nu}
+ \CC_{\delta\mu\psi\gamma\omega}A_{\psi \beta\mu}A_{\omega \delta\nu}\\
&+ \CC_{\omega\beta\psi\delta}A_{\psi \gamma\mu}A_{\omega \alpha\nu}
+ \CC_{\alpha\omega\psi\delta}A_{\psi \gamma\mu}A_{\omega \beta\nu}
+ \CC_{\alpha\beta\psi\delta}B_{\psi \gamma\mu\nu}
+ \CC_{\alpha\beta\psi\omega}A_{\psi \gamma\mu}A_{\omega \delta\nu}\\
&+ \CC_{\omega\beta\gamma\psi}A_{\psi \delta\mu}A_{\omega \alpha\nu}
+ \CC_{\alpha\omega\gamma\psi}A_{\psi \delta\mu}A_{\omega \beta\nu}
+ \CC_{\alpha\beta\omega\psi}A_{\psi \delta\mu}A_{\omega \gamma\nu}
+ \CC_{\alpha\beta\gamma\psi}B_{\psi \delta\mu\nu}
\Big)\varepsilon_{\gamma\delta}\\
&+\left(\CC_{\alpha\beta\gamma\delta,\mu}
+ \CC_{\psi\beta\gamma\delta}A_{\psi \alpha\mu}
+ \CC_{\alpha\psi\gamma\delta}A_{\psi \beta\mu}
+ \CC_{\alpha\beta\psi\delta}A_{\psi \gamma\mu}
+ \CC_{\alpha\beta\gamma\psi}A_{\psi \delta\mu}
\right)\varepsilon_{\gamma\delta,\nu}\\
&+\left(\CC_{\alpha\beta\gamma\delta,\nu}
+ \CC_{\psi \beta\gamma\delta}A_{\psi \alpha\nu}
+ \CC_{\alpha\psi\gamma\delta}A_{\psi \beta\nu}
+ \CC_{\alpha\beta\psi\delta}A_{\psi \gamma\nu}
+ \CC_{\alpha\beta\gamma\psi}A_{\psi \delta\nu}
\right)\varepsilon_{\gamma\delta,\mu}\\
&+ \CC_{\alpha\beta\gamma\delta} \varepsilon_{\gamma\delta,\mu\nu}\,,\label{eq:der2stress_local_comp}
\end{split}
\end{align}
\end{subequations}
where $\varepsilon_{\gamma\delta}$, $\varepsilon_{\gamma\delta,\mu}$, and  $\varepsilon_{\gamma\delta,\mu\nu}$ are
defined in \eqref{eq:strain_derstrain_local_comp}, while $A_{\psi \alpha\mu}$ and $B_{\psi \alpha\mu\nu}$ are:
\begin{subequations}\label{eq:A_B_local}
\begin{align}
A_{\psi \alpha\mu} &= \cfrac{\diff\a_\psi}{\diff x_\mu}\cdot\e_\alpha\,,\label{eq:A_local}\\
B_{\psi \alpha\mu\nu} &= \cfrac{\diff^2\a_\psi}{\diff x_\mu \diff x_\nu}\cdot\e_\alpha\,.\label{eq:B_local}
\end{align}
\end{subequations}
In the case that the basis $\{\a_\alpha\}$ depends on the parametric coordinates $(\xi^1,\xi^2,\xi^3)$, its derivatives with respect to the coordinates $\{x_\alpha\}$ can be computed by applying the chain rule. Thus, we first evaluate $\diff\a_\alpha/\diff\xi^\theta$ and  $\diff^2\a_\alpha/\diff\xi^\theta\diff\xi^\phi$ terms,
and from that, the quantities $\diff\a_\alpha/\diff x_\mu$ and  $\diff^2\a_\alpha/\diff x_\mu\diff x_\nu$
are obtained as herein detailed: \begin{subequations}
\begin{align}
\frac{\diff\a_\alpha}{\diff\xi^\theta} &= \frac{\diff\a_\alpha}{\diff x_\mu}\frac{\diff x_\mu}{\diff X_i}\frac{\diff X_i}{\diff \xi^\theta}\,,\label{eq:dedxi_local}\\
\begin{split}
\frac{\diff^2\a_\alpha}{\diff\xi^\theta\diff\xi^\phi} &=
\left(
 \frac{\diff^2\a_\alpha}{\diff x_\mu\diff x_\nu}\frac{\diff x_\mu}{\diff X_i}\frac{\diff x_\nu}{\diff X_j}
+\frac{\diff\a_\alpha}{\diff x_\mu}\frac{\diff^2 x_\mu}{\diff X_i\diff X_j}
\right)
 \frac{\diff X_i}{\diff \xi^\theta}\frac{\diff X_j}{\diff \xi^\phi}\\
&+ \frac{\diff\a_\alpha}{\diff x_\rho}\frac{\diff x_\rho}{\diff X_k}\frac{\diff^2 X_k}{\diff \xi^\theta\diff\xi^\phi}\,.
\end{split}\label{eq:de2dxi2_local}
\end{align}\label{eq:dedxi_de2dxi2_local}
\end{subequations}
Therefore, the derivatives $\diff\a_\alpha/\diff x_\mu$ and  $\diff^2\a_\alpha/\diff x_\mu\diff x_\nu$ are computed as:
\begin{subequations} \label{eq:dedx_d2edx2_local}
\begin{align}
\frac{\diff\a_\alpha}{\diff x_\mu} &= \frac{\diff\a_\alpha}{\diff \xi^\theta}\frac{\diff \xi^\theta}{\diff X_i}C_{i\mu}\,,\label{eq:dedx_local}\\
\frac{\diff^2\a_\alpha}{\diff x_\mu x_\nu} &=
\left( \frac{\diff^2\a_\alpha}{\diff\xi^\theta\diff\xi^\phi}
- \frac{\diff\a_\alpha}{\diff x_\rho}C_{k\rho}\frac{\diff^2 X_k}{\diff\xi^\theta\xi^\phi}\right)
\frac{\diff\xi^\theta}{\diff X_i}\frac{\diff\xi^\phi}{\diff X_j} C_{i\mu}C_{j\nu}\,,\label{eq:d2edx2_local}
\end{align}
\end{subequations}
where the term $\diff^2 x_\mu/\diff X_i\diff X_j$, according to Eqs.~\eqref{eq:metric} and \eqref{eq:metric_ders}, vanished and $\diff \xi^\theta /\diff X_i$ correspond to the inverse of the derivative of the geometric mapping, i.e.:
\begin{align}\label{eq:dXdxi_comp}
\frac{\diff X_i}{\diff \xi^\theta} = \frac{\diff\bm{F}}{\diff\xi^\theta}\cdot\E_i\,,
\end{align}
where $\bm{F}$ is the geometrical mapping that defines the physical domain $\Omega$ (previously introduced in Section~\ref{subsec:kin}).
Thus, the vectors $\{\diff \bm{F} / \diff \xi^1,\, \diff \bm{F} / \diff \xi^2,\,\diff \bm{F} / \diff \xi^3\}$ constitute the covariant basis of $\bm{F}$.
In the same way, $\diff X_i /\diff \xi^\theta\diff \xi^\phi$ corresponds to its second derivative:
\begin{align}\label{eq:d2Xdxi2_comp}
\frac{\diff^2 X_i}{\diff\xi^\theta\diff\xi^\phi} = \frac{\diff^2\bm{F}}{\diff\xi^\theta\diff\xi^\phi}\cdot\E_i\,.
\end{align}

Finally, the derivatives $\diff\a_\alpha/\diff \xi^\theta$ and $\diff^2\a_\alpha/\diff \xi^\theta\diff \xi^\phi$, that appear in Eq.~\eqref{eq:dedx_d2edx2_local}, are detailed in Appendix~\ref{subsec:localBasis_der} for a particular choice of the basis $\{\a_\alpha\}$.

\subsection{Stress divergence respect to the global reference system}\label{subsec:divergence_global_ref}
To fully address the system of Equations \eqref{eq:strong_form} in the global reference system $\{X_i\}$, we detail the stress tensor divergence components  \eqref{eq:stress_div_global_comp} as
\begin{equation}\label{eq:stress_div_global_comp_detail}
\begin{aligned}
\tilde\sigma_{ij,j} =&\frac{\diff \tilde\sigma_{ij}}{\diff X_j}=
\Big(D_{i\alpha}\,D_{j\beta}\,D_{k\gamma}\,D_{l\delta}\,\frac{\diff \CC_{\alpha\beta\gamma\delta}}{\diff X_j}+
\tilde A_{i \alpha j}\,D_{j\beta}\,D_{k\gamma}\,D_{l\delta}\,\CC_{\alpha\beta\gamma\delta}+D_{i\alpha}\,\tilde A_{j \beta j}\,D_{k\gamma}\,D_{l\delta}\,\CC_{\alpha\beta\gamma\delta}\\
&+D_{i\alpha}\,D_{j\beta}\,\tilde A_{k \gamma j}\,D_{l\delta}\,\CC_{\alpha\beta\gamma\delta}+D_{i\alpha}\,D_{j\beta}\,D_{k\gamma}\,\tilde A_{l \delta j}\,\CC_{\alpha\beta\gamma\delta}\Big)\tilde\varepsilon_{kl}+\tilde{\CC}_{ijkl}\,\tilde\varepsilon_{kl,j}\,,
\end{aligned}
\end{equation}
where $\tilde A_{i \alpha j}$ is defined as:
\begin{align}
\tilde A_{i \alpha j}=\dfrac{\diff D_{i\alpha}}{\diff X_j} = \dfrac{\diff\a_\alpha}{\diff X_j}\cdot\E_{i} = \dfrac{\diff\a_\alpha}{\diff x_\mu}\cdot\E_{i}\,C_{j\mu}\,.
\end{align}
The computation of $\diff\a_\alpha/\diff x_\mu$ was already defined in Eq.~\eqref{eq:dedx_local}.
Finally, applying the chain rule, the term $\diff \CC_{\alpha\beta\gamma\delta} / \diff X_j$ can be computed from $\diff \CC_{\alpha\beta\gamma\delta}/\diff x_\mu$ as:
\begin{align}
\frac{\diff \CC_{\alpha\beta\gamma\delta}}{\diff X_j} = \frac{\diff \CC_{\alpha\beta\gamma\delta}}{\diff x_\mu}\,C_{j\mu}\,.
\end{align}

	\section{A pointwise local basis for anisotropic materials}\label{subsec:localBasis_der}
In order to define the fibers orientation for anisotropic materials (including orthotropic ones) an appropriate local basis $\{\a_\alpha\}$ needs to be set.
A particular choice for this basis is an ortho-normalized version of the covariant basis such that:
\begin{subequations}\label{eq:localBasis}
\begin{align}
  \a_1 &= \frac{\g_1}{\lVert\g_1\rVert}\,,\label{eq:localBasis_e1}\\
  \a_2 &= \a_3\times\a_1\,,\label{eq:localBasis_e2}\\
  \a_3 &= \frac{\g_1\times\g_2}{\lVert\g_1\times\g_2\rVert}\,,\label{eq:localBasis_e3}
\end{align}
\end{subequations}
where $\{\g_1,\g_2,\g_3\}$ is the covariant basis, i.e.:
\begin{align}\label{eq:cov_vecs}
\g_\theta=\frac{\diff\bm{F}}{\diff\xi^\theta}\,.
\end{align}
Then, the first derivatives of the basis $\{\a_\alpha\}$ with respect to the parametric coordinates are:
\begin{subequations}\label{eq:dedxi_eval}
\begin{align}
  \frac{\diff\a_1}{\diff\xi^\theta} &= \frac{1}{\lVert\g_1\rVert}\left(\one-\a_1\otimes\a_1\right)\frac{\diff \g_1}{\diff\xi^\theta}\,,\label{eq:dedxi_eval_1}\\
  \frac{\diff\a_2}{\diff\xi^\theta} &= \frac{\diff\a_3}{\diff\xi^\theta}\times\a_1 + \a_3\times\frac{\diff\a_1}{\diff\xi^\theta}\,,\label{eq:dedxi_eval_2}\\
  \frac{\diff\a_3}{\diff\xi^\theta} &= \frac{1}{\lVert\g_1\times\g_2\rVert}\left(\one-\a_3\otimes\a_3\right)\left(\frac{\diff\g_1}{\diff\xi^\theta}\times\g_2 + \g_1\times\frac{\diff\g_2}{\diff\xi^\theta}\right)\,,\label{eq:dedxi_eval_3}
\end{align}
\end{subequations}
where $\bm{I}$ is the identity tensor.
Finally, the second derivatives are computed as:
\begin{subequations}\label{eq:d2edxi2_eval}
\begin{align}
\begin{split}\label{eq:d2edxi2_eval_1}
  \frac{\diff^2\a_1}{\diff\xi^\theta\diff\xi^\phi} &= \frac{1}{\lVert\g_1\rVert}
	\Bigg[-\left(\frac{\diff\a_1}{\diff\xi^\phi}\otimes\a_1+\a_1\otimes\frac{\diff\a_1}{\diff\xi^\phi}\right)\frac{\diff\g_1}{\diff\xi^\theta}\\
	&+\left(\one-\a_1\otimes\a_1\right)\frac{\diff^2\g_1}{\diff\xi^\theta\diff\xi^\phi}
	-\left(\a_1\cdot\frac{\diff\g_1}{\diff\xi^\phi}\right)\frac{\diff\a_1}{\diff\xi^\theta}\Bigg]\,,\\
\end{split}\\
\begin{split}\label{eq:d2edxi2_eval_2}
  \frac{\diff^2\a_2}{\diff\xi^\theta\diff\xi^\phi} &=
	\frac{\diff^2\a_3}{\diff\xi^\theta\diff\xi^\phi}\times\a_1
	+\frac{\diff\a_3}{\diff\xi^\theta}\times\frac{\diff\a_1}{\diff\xi^\phi}
	+\frac{\diff\a_3}{\diff\xi^\phi}\times\frac{\diff\a_1}{\diff\xi^\theta}
	+\a_3\times\frac{\diff^2\a_1}{\diff\xi^\theta\diff\xi^\phi}\,,\\
\end{split}\\
\begin{split}\label{eq:d2edxi2_eval_3}
  \frac{\diff^2\a_3}{\diff\xi^\theta\diff\xi^\phi} &= \frac{1}{\lVert\g_1\times\g_2\rVert}
	\Bigg[-\left(\frac{\diff\a_3}{\diff\xi^\phi}\otimes\a_3+\a_3\otimes\frac{\diff\a_3}{\diff\xi^\phi}\right)
	\left(\frac{\diff\g_1}{\diff\xi^\theta}\times\g_2+\g_1\times\frac{\diff\g_2}{\diff\xi^\theta}\right)\\
	&+\left(\one-\a_3\otimes\a_3\right)
	\left(\frac{\diff^2\g_1}{\diff\xi^\theta\diff\xi^\phi}\times\g_2
	+ \frac{\diff\g_1}{\diff\xi^\theta}\times\frac{\diff\g_2}{\diff\xi^\phi}
	+ \frac{\diff\g_1}{\diff\xi^\phi}\times\frac{\diff\g_2}{\diff\xi^\theta}
	+ \g_1\times\frac{\diff^2\g_2}{\diff\xi^\theta\diff\xi^\phi}\right)\\
	&-\left(\a_3\cdot\left(\frac{\diff\g_1}{\diff\xi^\phi}\times\g_2+\g_1\times\frac{\diff\g_2}{\diff\xi^\phi}\right)\right)\frac{\diff\a_3}{\diff\xi^\theta}\Bigg]\,.
\end{split}
\end{align}
\end{subequations}

    \bibliographystyle{unsrt}
    \bibliography{references}	
    
\end{document}